\documentclass[11pt,a4paper]{article}

\usepackage{amsmath,amssymb}
\newcommand{\e}{\varepsilon}
\newcommand{\va}{\varphi}
\newcommand{\D}{\Delta}
\newcommand{\La}{\Lambda}
\newcommand{\n}{\nabla}

\newcommand{\N}{\frac{N}{2}}
\newcommand{\NN}{\frac{N}{p}}

\newcommand{\p}{\partial}

\newcommand{\R}{\mathbb{R}}
\newcommand{\h}{\hookrightarrow}

\newcommand{\del}{\bar{\delta}}
\newcommand{\ka}{\bar{\kappa}}

\def\tilde{\widetilde}
\def\hat{\widehat}

\newtheorem{definition}{Definition}
\newtheorem{theorem}{Theorem}
\newtheorem{proposition}{Proposition}

\newtheorem{corollaire}{Corollary}

\newtheorem{remarka}{Remark}

\title{Strong solution for Korteweg system in $\mbox{bmo}^{-1}(\R^N)$ with initial density in $L^\infty$}
\author{Boris Haspot  \thanks{Universit\'e Paris Dauphine, PSL Research University, Ceremade, Umr Cnrs 7534, Place du Mar\' echal De Lattre De Tassigny 75775 Paris cedex 16 (France), haspot@ceremade.dauphine.fr } \thanks{ANGE project-team (Inria, Cerema, UPMC, CNRS), %
2 rue Simone Iff, CS 42112, 75589  Paris, France.}}
\date{}
\begin{document}
\maketitle
\begin{abstract}
In this paper we investigate the question of the local existence of strong solution for the Korteweg system in critical spaces in dimension $N\geq 1$ provided that the initial data are small. More precisely the initial momentum $\rho_0 u_0$ belongs to $\mbox{bmo}_{T}^{-1}(\R^N)$ for $T>0$ and the initial density $\rho_0$ is in $L^\infty(\R^N)$ and far away from the vacuum. This result extends the so called Koch-Tataru Theorem for the incompressible Navier-Stokes equations to the case of the Korteweg system. It is also interesting to observe that any initial shock on the density is instantaneously regularized inasmuch as the density becomes Lipschitz for any $\rho(t,\cdot)$ with $t>0$. We also prove the existence of global strong solution for small initial data $(\rho_0-1,\rho_0u_0)$ in the homogeneous Besov spaces $(\dot{B}^{\N-1}_{2,\infty}(\R^N)\cap \dot{B}^{\N}_{2,\infty}(\R^N)\cap L^\infty(\R^N))\times (\dot{B}^{\N-1}_{2,\infty}(\R^N))^N$. This result allows in particular to extend in dimension $N=2$ the notion of Oseen solutions defined for incompressible Navier-Stokes equations to the case of the Korteweg system when the vorticity of the momentum $\rho_0 u_0$ is a Dirac mass $\alpha\delta_0$ with $\alpha$ sufficiently small.
\end{abstract}
\section{Introduction}
We are concerned with compressible fluids endowed with internal
capillarity. The model we consider  originates from the XIXth
century work by J. F. Van der Waals and D. J. Korteweg \cite{VW,fK} and was
actually derived in its modern form in the 1980s using the second
gradient theory (see \cite{fDS,fJL,fTN}). 
Korteweg-type models are based on an extended version of
nonequilibrium thermodynamics, which assumes that the energy of the
fluid not only depends on standard variables but also on the
gradient of
the density. 
\\
We are now going to consider the so-called local Korteweg system which is a
compressible capillary fluid model, it can be derived from a Cahn-Hilliard like free energy (see the
pioneering work by J.- E. Dunn and J. Serrin in \cite{fDS} and also 
\cite{fA,fC,fGP}).
The conservation of mass and of momentum write:
\begin{equation}
\begin{cases}
\begin{aligned}
&\frac{\p}{\p t}\rho+{\rm div}(\rho u)=0,\\
&\frac{\p}{\p t}(\rho u)+{\rm div}(\rho
u\otimes u)-\rm div(2\mu \rho \,D (u))+\n P(\rho)={\rm div}K.
\end{aligned}
\end{cases}
\label{3systeme}
\end{equation}
Here $u=u(t,x)\in\R^{N}$ stands for the velocity field with $N\geq 1$, $\rho=\rho(t,x)\in\R^{+}$ is the density, $D(u)=\frac{1}{2}(\n u+^{t}\n u)$ is the strain tensor and $P(\rho)$ is the pressure (we will only consider regular pressure law). We denote by $\mu>0$ the viscosity coefficients of the fluid. We supplement the problem with initial condition $(\rho_{0},u_{0})$. The Korteweg tensor reads as:
\begin{equation}
{\rm div}K
=2 \kappa^2 \rho\n(\frac{\D\sqrt{\rho}}{\sqrt{\rho}}),
\label{divK}
\end{equation}
where $\kappa(\rho)=\frac{\kappa^2}{\rho}$ is the coefficient of capillarity with $\kappa^2 \in\R^+$. The term
${\rm div}K$  allows to describe the variation of density at the interfaces between two phases, generally a mixture liquid-vapor. In our case the capillary term is also called quantum pressure.\\[2mm]
We briefly mention that the existence of global strong solutions for the system (\ref{3systeme}) with small initial data for $N\geq2$ is known since the works by Hattori and Li \cite{Hat} in the case of constant capillary coefficient $\kappa(\rho)$. Danchin and Desjardins in \cite{fDD} improved this result by working with initial data $(\rho_{0}-1,\rho_{0}u_{0})$
belonging to the homogeneous Besov spaces \footnote{In the sequel we will only consider homogeneous Besov space that we will note $B^{s}_{p,r}$ even if they are generally written $\dot{B}^{s}_{p,r}$ with $s\in\R$ and $(p,r)\in[1,+\infty]^2$. We refer to \cite{BCD} for the definition of the homogeneous Besov spaces and for the notion of product in the Besov spaces related to the so called paraproduct. We will use the same notation as in \cite{BCD}.} $B^{\frac{N}{2}}_{2,1}\times B^{\frac{N}{2}-1}_{2,1}$ (it is important to point out that $B^{\frac{N}{2}}_{2,1}$ is embedded in $L^{\infty}$ which allows to control the $L^{\infty}$ norm of the density, it is even better since it implies that $\rho_0$ is necessary a continuous function).
In \cite{JMAA}, it is proved the existence of global strong solution with small initial data provided that $(\rho_{0}-1,\rho_{0}u_{0})$
is in $(B^{\frac{N}{2}}_{2,2}\cap L^\infty)\times B^{\frac{N}{2}-1}_{2,2}$. This result extends \cite{fDD} but does not allow to deal with general shocks on the initial density.
\\[2mm]
We wish now rewrite the system (\ref{3systeme}) using the formulation introduced in \cite{Aa,Ju} (where the existence of global energy weak solution is proved), we consider then the following effective velocity with $c>0$:
\begin{equation}
v_c=u+c\n\ln\rho,
\end{equation}
which enables us to rewrite the system (\ref{3systeme}) as follows when $\kappa^2\leq \mu^2$  (we will only consider this case in the sequel) using the fact that ${\rm div}(\rho\n\n\ln\rho)=2\rho\n (\frac{\D\sqrt{\rho}}{\sqrt{\rho}})$:
\begin{equation}
\begin{cases}
\begin{aligned}
&\p_t\rho-c\D\rho+{\rm div}(\rho v_c)=0\\
&\rho\p_t v_c+\rho u\cdot\n v_c-\mu{\rm div}(\rho\n v_c)-(\mu-c){\rm div}(\rho ^t\n v_c)+\n P(\rho)-2\kappa_1^2 \rho\n(\frac{\D\sqrt{\rho}}{\sqrt{\rho}})=0.
\end{aligned}
\label{k2}
\end{cases}
\end{equation}
with $\kappa_1^2=\kappa^2-2\mu c +c^2$. 
We specify now the value of $c$ and we want to deal with a $c$ such that $\kappa^2_1=0$, we take then:
$$c_1=\mu-\sqrt{\mu^2-\kappa^2}\;\;\mbox{and}\;\;c_2=\mu+\sqrt{\mu^2-\kappa^2}.$$
Setting now $v_1=u+c_1\n\ln\rho$ and $v_2=u+c_2\n\ln\rho$ , we have:
\begin{equation}
\begin{cases}
\begin{aligned}
&\p_t\rho-c_1\D\rho+{\rm div}(\rho v_1)=0\\
&\p_t(\rho v_1)+\frac{1}{2}{\rm div}(\rho v_1\otimes v_2)+\frac{1}{2}{\rm div}(\rho v_2\otimes v_1)-\mu\D (\rho v_1)\\
&\hspace{5cm}-\sqrt{\mu^2-\kappa^2}\n{\rm div}(\rho v_1)+\n P(\rho)=0,\\[2mm]
&\p_t(\rho v_2)+\frac{1}{2}{\rm div}(\rho v_1\otimes v_2)+\frac{1}{2}{\rm div}(\rho v_2\otimes v_1)-\mu\D (\rho v_2)\\
&\hspace{5cm}+\sqrt{\mu^2-\kappa^2}\n{\rm div}(\rho v_2)+\n P(\rho)=0.
\end{aligned}
\label{kfinal}
\end{cases}
\end{equation}
In this paper our main goal consists in proving the existence of global or local strong solutions for the system (\ref{3systeme}) with minimal regularity assumption on the initial data. More precisely since this system models a mixture liquid vapor with different density, we wish to show the existence of strong solution for initial density $\rho_0$ belonging only to the set $L^\infty(\R^N)$. In particular it implies that our functional setting will include the case of initial density admitting shocks what is also fundamental both in the physical theory of non-equilibrium thermodynamics as well as in the mathematical study of inviscid models for compressible flow. In addition we wish also to deal with momentum $\rho_0 u_0$ exhibiting specific structure, typically we have in mind the case of initial vorticity belonging to the set of finite measure in dimension $N=2$ (we will recall later that there is such a theory for Navier-Stokes equations and even explicit such solutions, the so called Oseen tourbillon). In order to obtain such results, it seems necessary to work in space with minimal regularity assumptions. Typically the space  $B^{\NN-1}_{p,r}$ with third index $r=+\infty$ is a good candidate for the initial momentum $\rho_0u_0$.
\\
To do this, let us now recall the notion of scaling for the Korteweg's system (\ref{3systeme}). Such an
approach is now classical for incompressible Navier-Stokes equations
and yields local well-posedness (or global well-posedness for small
initial data) in spaces with minimal regularity.
In our situation we can
easily check that, if $(\rho,u)$ solves (\ref{3systeme}), then
$(\rho_{\lambda},u_{\lambda})$ solves also this system:
\begin{equation}
\rho_{\lambda}(t,x)=\rho(\lambda^{2}t,\lambda x)\,,\,u_{\lambda}(t,x)=\lambda u(\lambda^{2}t,\lambda x)\;\;\;\mbox{with}\;\lambda\in\R^+
\label{scal1}
\end{equation}
provided the pressure laws $P$ have been changed to
$\lambda^{2}P$.
\begin{remarka}
It is very important to point out that since there is only a scaling invariance up to the pressure term, we can not hope to show the existence of global strong self similar solution for the Korteweg system. Indeed in comparison with the Navier-Stokes system, we know that it exists global in time self similar solution provided that the initial velocity $u_0$ is homogeneous of degree $-1$ (such initial data exists for example in the set $B^{\NN-1}_{p,\infty}$ with $1\leq p<+\infty$). In particular in the sequel we will able to deal with initial data $(\rho_0,\rho_0 u_0)$ such that $\rho_0$ is homogeneous of degree $0$ and $\rho_0 u_0$ of degree $-1$, however the associated strong solution will be a priori not self similar.
\end{remarka}
The previous transformation (\ref{scal1}) suggests us however to choose initial data
$(\rho_{0},u_{0})$ in spaces whose norm is invariant for
all $\lambda>0$ by the transformation
$(\rho_{0},u_{0})\longrightarrow(\rho_{0}(\lambda\cdot),\lambda
u_{0}(\lambda\cdot)).$ 
A natural candidate is for example the space $L^\infty(\R^N)\times \dot{H}^{\N-1}$. As it was mentioned, this invariance was also used initially by Kato \cite{K1} to prove that the Navier-Stokes system is locally well-posed for arbitrary data in $L^N( \R^N)$ (when $N\geq 2$) and  globally well-posed for  small initial data. Kato's result was extended to larger scale invariant function spaces (one interest of dealing with larger function spaces is that they may contain initial data which are homogeneous of degree $-1$ and therefore give rise to self-similar solutions). In particular in \cite{CMP} Cannone, Meyer and Planchon proved the existence of global strong solution with small initial data in $B^{\NN-1}_{p,\infty}$ with $p<+\infty$. A similar analysis was carried out for the vorticity equation in Morrey spaces by Giga and Miyakawa \cite{GM}.
 Finally this approach has been generalized by 
 Koch and Tataru in \cite{KT} when the initial data is small in $BMO^{-1}(\R^N)$. \\
These results allow in particular to obtain the existence of global self similar solution for small initial data when $N\geq 2$ for the incompressible Navier-Stokes equations. In dimension $N=2$, it proves the existence and the uniqueness of solution for initial data $u_0$ satisfying ${\rm curl}u_0=\alpha\delta_0$ (which is equivalent to $u_0=\alpha\frac{x^{\perp}}{|x|^2}$) with $\alpha$ small enough. These solutions are the so-called Lamb-Oseen solutions which are self-similar. Let us emphasize in addition that we have even an explicit formula for these solutions even when $|\alpha|$ is large, the Lamb-Oseen vortex are given by:
\begin{equation}
\begin{aligned}
&{\rm curl}u_\alpha (t,x)=\frac{\alpha}{t}G(\frac{x}{\sqrt{t}}),\;u_\alpha(t,x)=\frac{\alpha}{\sqrt{t}}v_{G}(\frac{x}{\sqrt{t}}),\;x\in\R^2,\;t>0,
\end{aligned}
\label{Oseen}
\end{equation}
where:
$$G(\xi)=\frac{1}{4\pi}e^{-\frac{|\xi|^2}{4}},\;v_G(\xi)=\frac{1}{2\pi}\frac{\xi^{\perp}}{|\xi|^2}(1-e^{-\frac{|\xi|^2}{4}}),\;\xi\in\R^2,$$
with $\xi^{\perp}=(-\xi_2,\xi_1)$. In passing, we mention that the existence of global weak solution with  initial vorticity in the set ${\cal M}(\R^2)$ of all finite real measures on $\R^2$ was first proved by Cottet \cite{C} and independently by Giga, Miyakawa and Osada \cite{GM}. In \cite{GM}, the authors proved also the uniqueness when the atomic part of the initial vorticity is sufficiently small.
The uniqueness for any ${\rm curl}u_0\in{\cal M}(\R^2)$ is proved in \cite{GG}, it allows in particular  to obtain the existence and the uniqueness of global self similar solution for large initial data when $N=2$.\\
In this paper we are interested in extending the technics of Koch-Tataru to the case of the Korteweg system (\ref{kfinal}). More precisely, we wish to prove the existence of strong solution in finite time for (\ref{3systeme}) provided that $\rho_0 v_1(0,\cdot)$ and $\rho_0 v_2(0,\cdot)$ are sufficiently small in norm $\|\cdot\|_{bmo_T^{-1}(\R^N)}$ for $T>0$. Furthermore we will assume that the initial density is far away from the vacuum and is bounded in $L^\infty$ norm . We recall that
 $bmo^{-1}(\R^N)$ is the set of temperated distribution $u_0$ for which for all $T>0$ we have:
$$\sup_{x\in\R^N,0<t<T}\big( t^{-\N}\int^t_0\int_{B(x,\sqrt{t})}|e^{s\D}u_0 (y)|^2 dy ds\big)^{\frac{1}{2}}<+\infty.$$
We define the norm $\|\cdot\|_{bmo_T^{-1}(\R^N)}$ by:
$$\|u_0\|_{bmo_T^{-1}(\R^N)}=\sup_{x\in\R^N,0<t<T}\big( t^{-\N}\int^t_0\int_{B(x,\sqrt{t})}|e^{s\D}u_0 (y)|^2 dy ds\big)^{\frac{1}{2}},$$
with $e^{t\D}u_0$ be the solution to the heat equation with initial data $u_0$:
$$e^{t\D}u_0=u_0*\phi_{\sqrt{4t}}\;\;\mbox{with}\;\phi(x)=\pi^{-\N}e^{-|x|^2}\;\;\mbox{and}\;\;\phi_t (x)=t^{-N}\phi(\frac{x}{t}).$$
In the sequel we will denote by ${\cal E}_T$ the space of temperated distribution associated to the following norm:
\begin{equation}
\|u\|_{{\cal E}_T}=\sup_{0<t\leq T}\sqrt{t}\|u(t,\cdot)\|_{L^\infty(\R^N)}+\big(\sup_{x\in\R^N, 0<t\leq T}t^{-\N}\int^t_0\int_{ B(x,\sqrt{t})}|u|^2(s,y) dt dy\big)^{\frac{1}{2}},
\label{defXs}
\end{equation}
In the sequel, we will show that $\rho v_1$ and $\rho v_2$ verifying (\ref{kfinal}) belong to ${\cal E}_T$ for $T>0$ the time of existence and that $(\rho,\frac{1}{\rho})$ is in $L^\infty_T(L^\infty(\R^N))$. We recall that $\rho v_1=\rho u+c_1\n\rho$ and $\rho v_2=\rho u+c_2\n\rho$, it implies in particular that $\n\rho=\frac{1}{2\sqrt{\mu^2-\kappa^2}}(\rho v_2-\rho v_1)$ and we deduce from the definition of ${\cal E}_T$ that it exists $C>0$ independent on $t$ such that for all $0<t\leq T$ we have:
\begin{equation}
\|\n\rho(t,\cdot)\|_{L^\infty(\R^N)}\leq \frac{C}{\sqrt{t}}.
\label{regu}
\end{equation}
Combining (\ref{regu}) and the fact that $\rho$ remains in $L^\infty_T(L^\infty(\R^N))$ it implies that instantaneously the density $\rho$ is regularized and becomes necessary Lipschitz even if the initial density $\rho_0$ admits shocks.\\
The second point is that with our choice on the initial data, it implies that $\rho_0 u_0$ is in $bmo^{-1}_T$.  In other way we can work with Dirac mass $\alpha\delta_0$ for the initial vorticity associated to the momentum $\rho_0 u_0$ in dimension $N=2$ provided that $\alpha$ is sufficiently small, it consists to prescribe the initial momentum as follows $\rho_0 u_0(x)=\alpha\frac{x^\perp}{|x|^2}$ with $x\in\R^2$.
We can then extend the notion of Lam-Oseen tourbillon to the case of the Korteweg system at least when $\alpha$ is sufficiently small. It is obvious that in this compressible framework the divergence of such Lamb-Oseen tourbillon does not remain null all along the time, this is due to the coupling between vorticity and divergence of the velocity. Similarly we can deal with vorticity vortex and divergence vortex if we take the following initial momentum  $\rho_0 u_0(x)=\alpha\frac{x^\perp}{|x|^2}+\alpha_1\frac{x^\perp}{|x|^2} $ with $x\in\R^2$ and $|\alpha|$, $|\alpha_1|$ sufficiently small.
Finally we will extend the previous result in dimension $N\geq 2$ by proving the existence of global strong solution provided that the initial data $(\rho_0-1,\rho_0 u_0)$ are small enough in $(\dot{B}^{\N-1}_{2,\infty}(\R^N)\cap \dot{B}^{\N}_{2,\infty}(\R^N)\cap L^\infty(\R^N))\times (\dot{B}^{\N-1}_{2,\infty}(\R^N))^N$. It enables us again to obtain global strong Oseen solutions for the Korteweg system provided that the initial data is sufficiently small.
\subsection{Mathematical results}
In this section we state our main result. 
\begin{theorem}
Let $0<\kappa^2\leq \mu^2$ and $N\geq 1$.
Let $\rho_0\in L^\infty(\R^N)$ with $\rho_0\geq c>0$ and $\mathbb{P}(\rho_0 v_1(0,\cdot)), \mathbb{Q}(\rho_0 v_1(0,\cdot)), \mathbb{P}(\rho_0 v_2(0,\cdot)),  \mathbb{Q}(\rho_0 v_2(0,\cdot) )\in bmo^{-1}(\R^N)$ \footnote{ We define here $\mathbb{Q}=\n(\D)^{-1}{\rm div}$ and $\mathbb{P}=Id-\mathbb{P}$.}. Then there exists $T>0$ sufficiently small such that there exists a  unique solution $(\rho,\rho u)$ of the system (\ref{3systeme}) on $[0,T]$ provided that for $\e_1>0$ small enough we have :
\begin{equation}
\begin{aligned}
&\|\mathbb{P}(\rho_0 v_1(0,\cdot))\|_{bmo^{-1}_T}+\|\mathbb{Q}(\rho_0 v_1(0,\cdot))\|_{bmo^{-1}_T}+\|\mathbb{P}(\rho_0 v_2(0,\cdot))\|_{bmo^{-1}_T}\\
&\hspace{7,5cm}+\|\mathbb{Q}(\rho_0 v_2(0,\cdot))\|_{bmo^{-1}_T}\leq\e_1,
\end{aligned}
\label{smallness}
\end{equation}
In addition it exists $C>0$ such that:
 \begin{equation}
 \|\rho v_1\|_{{\cal E}_T}+ \|\rho v_2\|_{{\cal E}_T}+\|\rho\|_{L^\infty_T(L^\infty)}+\|\frac{1}{\rho}\|_{L^\infty_T(L^\infty)}\leq C.
 \label{estimation}
 \end{equation}
\label{theo1}
\end{theorem}
\begin{remarka}
It is important to observe that the condition (\ref{smallness}) implies that $\n\rho_0$ is small in $bmo^{-1}(\R^N)$ with the norm $\|\cdot\|_{bmo^{-1}_T}$ when $0<\kappa^2<\mu^2$. It is also interesting to observe that there is no smallness assumption on $\rho_0$ when $\kappa^2=\mu^2$ since in this case $\rho_0 v_1(0,\cdot)=\rho_0 v_2(0,\cdot)$. We deduce also that in this case $\rho_0 u_0$ is not necessary small in $bmo^{-1}$ but the smallness carries on the momentum of the effective velocity
$\rho_0v_1(0,\cdot)$ which describes the coupling between the velocity $u_0$ and the density $\rho_0$.
\end{remarka}
\begin{remarka}
It is important to mention that the initial data are defined by the momentum $\rho_0v_1(0,\cdot)$ and 
 $\rho_0v_2(0,\cdot)$, indeed the initial velocity $u_0$ is a priori not defined (however the momentum $m_0$ which is equal roughly speaking to $\rho_0u_0$ is well defined).
 \end{remarka}
\begin{corollaire}
\label{cor1}
Let $0<\kappa^2\leq \mu^2$ and $N\geq 1$. Let $\rho_0\in L^\infty(\R^N)$ with $\rho_0\geq c>0$ , assume that $\n \rho_0$ and $\rho_0 u_0$ are in $\overline{{\cal D}(\R^N)}^{bmo^{-1}}$ then there exists $T>0$ such that there exists a unique solution $(\rho,\rho u)$ of the system (\ref{3systeme}). We have in addition the estimate (\ref{estimation}).
\end{corollaire}
We prove now a result of global strong solution with small initial data.
\begin{theorem}
Let $0<\kappa^2\leq \mu^2$, $\bar{\rho}>0$, $N\geq 2$ and $P'(\bar{\rho})>0$. We assume that $(\rho_0-\bar{\rho},u_0)\in (B^{\N-1}_{2,\infty}\cap B^{\N}_{2,\infty}\cap L^\infty)\times (B^{\N-1}_{2,\infty})^N$. There exists $\e_0$ such that if:
\begin{equation}
\|\rho_0-\bar{\rho}\|_{B^{\N-1}_{2,\infty}\cap B^{\N}_{2,\infty}\cap L^\infty}+\|u_0\|_{B^{\N-1}_{2,\infty}}\leq\e_0,
\end{equation}
then there exists a global strong solution $(\rho,u)$ of the system (\ref{3systeme}). In addition it exists $C>0$ such that for $s_1\in(\frac{3}{4},1)$:
$$
\begin{cases}
\begin{aligned}
&\|(\rho v_1,\rho v_2)\|_{\widetilde{L}^\infty(\R^+,B^{\N-1}_{2,\infty})\cap \widetilde{L}^1(\R^+,B^{\N+1}_{2,\infty}) }+\|\rho-\bar{\rho}\|_{\widetilde{L}^\infty(\R^+,B^{\N-1}_{2,\infty}\cap B^{\N}_{2,\infty})\cap \widetilde{L}^1(\R^+,B^{\N+1}_{2,\infty}\cap B^{\N+2}_{2,\infty}) }\leq C,\\
&\|\rho-\bar{\rho}\|_{L^\infty(\R^+,L^\infty(\R^N))}\leq C\\
&\sup_{t\in\R^+} t^{\frac{s_1}{2}}\|(\rho(t,\cdot)-\bar{\rho},\rho v_1(t,\cdot),\rho v_2(t,\cdot)\|_{B^{\N-1+s_1}_{2,\infty}}\leq C.
\end{aligned}
\end{cases}
$$
\label{theo2}
\end{theorem}
\begin{remarka}
As previously, this theorem allows to prove the existence of global strong Oseen solution provided that ${\rm curl}(\rho_0 u_0)=\alpha\delta_0$ with $|\alpha|$ small enough in dimension $N=2$.\\
When $\kappa^2<\mu^2$ it would be possible to extend this result by working with Besov spaces constructed on $L^p$ Lebesgue spaces in high frequencies as in \cite{arma} for compressible Navier-Stokes equations.
\end{remarka}
\section{Proof of the Theorem \ref{theo1} and the Corollary \ref{cor1}}
We are going now to prove the Theorem \ref{theo1} in the case $0<\kappa^2<\mu^2$. The case $\kappa^2=\mu^2$ is similar except that $v=v_1=v_2$ and we apply the estimates in $bmo^{-1}_T$ on $\rho v$. From (\ref{kfinal}), we observe that we have for $t>0$:
$$
\begin{cases}
\begin{aligned}
&\rho(t,\cdot)=e^{c_1t\D}\rho_0-\int^t_0 e^{c_1(t-s)\D}{\rm div}(\rho v_1) ds,\\
&\mathbb{P}\rho v_1 (t,\cdot)=e^{\mu t\D}\mathbb{P}(\rho v_1(0,\cdot)) -\frac{1}{2}\int^t_0 e^{\mu(t-s)\D} \mathbb{P}({\rm div}(\rho v_1\otimes v_2)+{\rm div}(\rho v_2\otimes v_1))  ds,\\
&\mathbb{P}\rho v_2 (t,\cdot)=e^{\mu t\D}\mathbb{P}(\rho v_2(0,\cdot)) -\frac{1}{2}\int^t_0 e^{\mu (t-s)\D} \mathbb{P}({\rm div}(\rho v_1\otimes v_2)+{\rm div}(\rho v_2\otimes v_1))  ds,\\
&\mathbb{Q}\rho v_1 (t,\cdot)=e^{c_2 t\D}\mathbb{Q}(\rho v_1(0,\cdot)) -\int^t_0 e^{c_2(t-s)\D} \mathbb{Q}(F(\rho,v_1,v_2)(s))  ds,\\
&\mathbb{Q}\rho v_2 (t,\cdot)=e^{c_1 t\D}\mathbb{Q}(\rho v_2(0,\cdot)) -\int^t_0 e^{c_1(t-s)\D} \mathbb{Q}(F(\rho,v_1,v_2)(s))  ds,\\
\end{aligned}
\end{cases}
$$
with:
$$F(\rho,v_1,v_2)(s)=\frac{1}{2}({\rm div}(\rho v_1\otimes v_2)+{\rm div}(\rho v_2\otimes v_1))(s,\cdot)+\n P(\rho(s,\cdot)).$$
We are going now to work in the following space:
$$\|(\rho,v_1,v_2)\|_{X_T}=\|\rho v_1\|_{{\cal E}_T}+\|\rho v_2\|_{{\cal E}_T}+\|(\rho,\frac{1}{\rho})\|_{L^\infty(\R^N)}, $$
for $T$ small enough.
We observe in particular that for any $\rho\in L^\infty_T(L^\infty(\R^N))$ we have:
\begin{equation}
\|\rho m\|_{{\cal E}_T}\leq\|\rho\|_{L^\infty_T(L^\infty(\R^N))}\|m\|_{{\cal E}_T}.
\label{3a}
\end{equation}
We shall use in the sequel a contracting mapping argument with the function $\psi$ defined as follows:
\begin{equation}
\begin{aligned}
&\psi(\rho,\rho v_1,\rho v_2)=\left(\begin{array}{c}
e^{c_1t\D}\rho_0\\
e^{\mu t\D}\mathbb{P}(\rho v_1(0,\cdot))+e^{c_2 t\D}\mathbb{Q}(\rho v_1(0,\cdot))\\
e^{\mu t\D}\mathbb{P}(\rho v_2(0,\cdot))+e^{c_1 t\D}\mathbb{Q}(\rho v_2(0,\cdot))\\
\end{array}
\right)\\
&-\int_{0}^{t}\left(\begin{array}{c}
 e^{c_1(t-s)\D}{\rm div}(\rho v_1) \\
 \frac{1}{2} e^{\mu(t-s)\D} \mathbb{P}({\rm div}(\rho v_1\otimes v_2)+{\rm div}(\rho v_2\otimes v_1)) +e^{c_2(t-s)\D} \mathbb{Q}(F(\rho,v_1,v_2))\\
 \frac{1}{2} e^{\mu(t-s)\D} \mathbb{P}({\rm div}(\rho v_1\otimes v_2)+{\rm div}(\rho v_2\otimes v_1)) +e^{c_1(t-s)\D} \mathbb{Q}(F(\rho,v_1,v_2)) \\
\end{array}
\right)\ ds.
\end{aligned}
\end{equation}
More precisely we want to prove that $\psi$ is a contractive map from $E_{R,M,T}$ to $E_{R,M,T}$ with $R,M>0$:
\begin{equation}
\begin{aligned}
&E_{R,M,T}=\{(\rho,\rho v_1,\rho v_2)\in L^\infty_T(L^\infty(\R^N))\times{\cal E}_T\times{\cal E}_T / \\
&\rho\geq \frac{1}{M}\;\;\mbox{on}\;[0,T]\times\R^N, \,\|\rho v_1\|_{{\cal E}_{T}}+\|\rho v_2\|_{{\cal E}_{T}}\leq 2R, \|\rho\|_{L^\infty_T(L^\infty(\R^N)}\leq 2\|\rho_0\|_{L^\infty} \}
\end{aligned}
\label{def}
\end{equation}
$E_{R,M,T}$ is endowed with the following norm:
\begin{equation}
\|(\rho,\rho v_1,\rho v_2)\|_{E_{R,M,T}}=\|(\rho v_1,\rho v_2)\|_{{\cal E}_T}+\frac{1}{\beta}\|\rho\|_{L^\infty_T(L^\infty)},
\end{equation}
with $\beta>0$ sufficiently large that we will defined later and $E_{R,M,T}$ is a Banach space.\\
Let $c>0$, we know that it exists $C>0$ (depending on $c$) such that for all $0<t\leq T$ we have (see \cite{Lemarie} p163):
\begin{equation}
\begin{cases}
\begin{aligned}
&\|\int^t_0 e^{c(t-s)\D} \mathbb{P}{\rm div}(v\otimes w)ds\|_{{\cal E}_T}\leq C\|v\|_{{\cal E}_T}\|w\|_{{\cal E}_T}\\
&\|\int^t_0 e^{c(t-s)\D} \mathbb{Q} {\rm div}(v\otimes w)ds\|_{{\cal E}_T}\leq C\|v\|_{{\cal E}_T}\|w\|_{{\cal E}_T}.
\end{aligned}
\end{cases}
\label{estim1}
\end{equation}
Similarly we have (see \cite{KT}) for $C,C_1>0$:
\begin{equation}
\begin{aligned}
&\|\int^t_0 e^{c(t-s)\D} \n P(\rho)(s,\cdot) ds\|_{{\cal E}_T} \\
&\leq C\big(\sup_{0<t\leq T}t \|P(\rho)(t,\cdot)\|_{L^\infty}+\sup_{0<t\leq T,x\in\R^N}\frac{1}{t^\N}\int^t_0\int_{B(x,\sqrt{t})}|P(\rho)(s,y)| ds dy\big)\\
&\leq C_1 T \|P(\rho)\|_{L^\infty_T(L^\infty(\R^N))}.
\end{aligned}
\label{estim2}
\end{equation}
We deduce from (\ref{3a}), (\ref{estim1}), (\ref{estim2}) and the definition of the $bmo^{-1}_T$ that it exists $C>0$ such that for any $T>0$ we have:
\begin{equation}
\begin{aligned}
&\|(\psi_2(\rho,\rho v_1,\rho v_2),\psi_3(\rho,\rho v_1,\rho v_2)) \|_{{\cal E}_T^2}\leq  C\big(\|\mathbb{P}(\rho_0 v_1(0,\cdot))\|_{bmo^{-1}_T}+\|\mathbb{Q}(\rho_0 v_1(0,\cdot))\|_{bmo^{-1}_T}\\
&+\|\mathbb{P}(\rho_0 v_2(0,\cdot))\|_{bmo^{-1}_T}+\|\mathbb{Q}(\rho_0 v_2(0,\cdot))\|_{bmo^{-1}_T}
+\|\frac{1}{\rho}\|_{L^\infty_T(L^\infty)}\|\rho v_1\|_{{\cal E}_T}\|\rho v_2\|_{{\cal E}_T}\\
&\hspace{9cm}+T \|P(\rho)\|_{L^\infty_T(L^\infty(\R^N))}\big).
\end{aligned}
\label{estim3}
\end{equation}
Let us estimate now the $L^\infty$ norm of $(\psi_1(\rho,\rho v_1,\rho v_2),\frac{1}{\psi_1(\rho,\rho v_1,\rho v_2)})$, we have then for any $t\in[0,T]$ and using the maximum principle:
\begin{equation}
\begin{cases}
\begin{aligned}
&
\|\psi_1(\rho,\rho v_1,\rho v_2)(t,\cdot)\|_{L^\infty}\leq \|\rho_0\|_{L^\infty}+\|\int^t_0 e^{c_1(t-s)\D} {\rm div}(\rho v_1) ds\|_{L^\infty}\\
&
\psi_1(\rho,\rho v_1,\rho v_2)(t,\cdot)\geq \min _{x\in\R^N}\rho_0(x)-\|\int^t_0 e^{c_1(t-s)\D} {\rm div}(\rho v_1) ds\|_{L^\infty}.
\end{aligned}
\end{cases}
\label{estim4}
\end{equation}
We have now using integration by parts and for $0<t\leq T$ it exists $C,C_1,C_2,C_3>0$ such that:
\begin{equation}
\begin{aligned}
&|\int^t_0 e^{c_1(t-s)\D} {\rm div}(\rho v_1) ds|\leq |\int^t_0\int_{\R^N}\frac{1}{(4c_1\pi(t-s))^{\N}} e^{-\frac{|x-y|^2}{4c_1(t-s)}}{\rm div}(\rho v_1)(s,y) dy ds|\\
&\leq C\int^t_0\int_{\R^N}\frac{1}{(4\pi c_1(t-s))^{\N+\frac{1}{2}}} \frac{|x-y|}{\sqrt{4 c_1(t-s)}}e^{-\frac{|x-y|^2}{4c_1(t-s)}} |\rho v_1(s,y)| dy ds|\\
&\leq C_1\int^t_0 \frac{1}{\sqrt{t-s}}\| \frac{1}{(4\pi(t-s))^{\N}} \frac{|\cdot|}{\sqrt{4(t-s)}}e^{-\frac{|\cdot|^2}{4(t-s)}}\|_{L^1} \|\rho v_1(s,\cdot)\|_{L^\infty}  ds\\
&\leq C_2\|\rho v_1\|_{{\cal E}_T}  \int^t_0 \frac{1}{\sqrt{t-s}}\frac{1}{\sqrt{s}} ds\leq C_3 \|\rho v_1\|_{{\cal E}_T} .
\end{aligned}
\label{estim5}
\end{equation}
We deduce now from (\ref{estim4}) and (\ref{estim5}) that for $C>0$ large enough we have for $t\in[0,T]$:
\begin{equation}
\begin{cases}
\begin{aligned}
&
\|\psi_1(\rho,\rho v_1,\rho v_2)(t,\cdot)\|_{L^\infty}\leq \|\rho_0\|_{L^\infty}+C \|\rho v_1\|_{{\cal E}_T} \\
&\psi_1(\rho,\rho v_1,\rho v_2)(t,\cdot)\geq \min_{x\in\R^N} \rho_0(x)-C \|\rho v_1\|_{{\cal E}_T} .
\end{aligned}
\end{cases}
\label{estim6}
\end{equation}
Let us prove now that for $M$, $T$, $R$ well chosen, we have:
\begin{equation}
\psi(E_{R,M,T})\subset E_{R,M,T}.
\label{inclu}
\end{equation}
Taking $\frac{1}{M}=\frac{c}{2}$ and  $R$ such that $2CR\leq \min(\frac{c}{2},\|\rho_0\|_{L^\infty(\R^N)})$ (with $C>0$ defined in  (\ref{estim6})) and where $\min_{x\in\R^N} \rho_0(x)\geq c>0$, we deduce from (\ref{estim6}) that for any $t\in[0,T]$ and for any $(\rho,\rho v_1,\rho v_2)\in E_{R,M,T}$:
\begin{equation}
\begin{cases}
\begin{aligned}
&
\|\psi_1(\rho,\rho v_1,\rho v_2)(t,\cdot)\|_{L^\infty}\leq 2 \|\rho_0\|_{L^\infty}\\
&\psi_1(\rho,\rho v_1,\rho v_2)(t,\cdot)\geq \frac{c}{2}=\frac{1}{M}.
\label{estim7}
\end{aligned}
\end{cases}
\end{equation}
Now we deduce from (\ref{estim3}) that we have  for any $(\rho,\rho v_1,\rho v_2)\in E_{R,M,T}$ and $C_1>0$ large enough:
\begin{equation}
\begin{aligned}
&\|(\psi_2(\rho,\rho v_1,\rho v_2),\psi_3(\rho,\rho v_1,\rho v_2)) \|_{{\cal E}_T^2}\leq  C_1\big( \mathbb{P}(\rho_0 v_1(0,\cdot))\|_{bmo^{-1}_T}+\|\mathbb{Q}(\rho_0 v_1(0,\cdot))\|_{bmo^{-1}_T}\\
&+\|\mathbb{P}(\rho_0 v_2(0,\cdot))\|_{bmo^{-1}_T}+\|\mathbb{Q}(\rho_0 v_2(0,\cdot))\|_{bmo^{-1}_T}
+4M R^2+TM_1\big),
\end{aligned}
\label{estim7a}
\end{equation}
with and $M_1=\sup_{x\in [0,2 \|\rho_0\|_{L^\infty}]} P(x)$. We set now: 
$$
\begin{aligned}
&R= 4C_1\big( \mathbb{P}(\rho_0 v_1(0,\cdot))\|_{bmo^{-1}_T}+\|\mathbb{Q}(\rho_0 v_1(0,\cdot))\|_{bmo^{-1}_T}+\|\mathbb{P}(\rho_0 v_2(0,\cdot))\|_{bmo^{-1}_T}\\
&\hspace{8cm}+\|\mathbb{Q}(\rho_0 v_2(0,\cdot))\|_{bmo^{-1}_T}\big).
\end{aligned}
$$
Now choosing $R$ such that $4C_1 MR^2\leq \frac{R}{2}$ and $T$ such that $C_1M_1 T \leq\frac{R}{4}$ we deduce from (\ref{estim7}) and (\ref{estim7a}) that we have  for any $(\rho,\rho v_1,\rho v_2)\in E_{R,M,T}$:
\begin{equation}
\begin{aligned}
&\|(\psi_2(\rho,\rho v_1,\rho v_2),\psi_3(\rho,\rho v_1,\rho v_2)) \|_{{\cal E}_T^2}\leq 2R\;\;\mbox{and}\;\;\psi(E_{R,M,T})\subset E_{R,M,T}..
\end{aligned}
\label{estim8}
\end{equation}
In conclusion we have chosen $R$ and $T$ such that:
\begin{equation}
\begin{aligned}
&R\leq \frac{\min(\frac{c}{2},\|\rho_0\|_{L^\infty(\R^N)})  }{2C},\;R\leq\frac{c}{16 C_1 }\;\;\;\mbox{and}\;\;\;T\leq\frac{R}{4C_1M_1  }.
\end{aligned}
\end{equation}
Let us prove now that $\psi$ is a contractive map. More precisely let us estimate $\|\psi(\rho,\rho v_1,\rho v_2)-\psi(\rho_1,\rho_1 w_1,\rho_1 w_2)\|_{E_{R,M,T}}$ with $((\rho,\rho v_1,\rho v_2),(\rho_1,\rho_1 w_1,\rho_1 w_2))\in E_{R,M,T}$. We have then:
\begin{equation}
\begin{aligned}
&\psi(\rho,\rho v_1,\rho v_2)-\psi(\rho_1,\rho_1 w_1,\rho_1 w_2)=\\
&\int_{0}^{t}\left(\begin{array}{c}
 e^{c_1(t-s)\D}{\rm div}(\rho_1w_1-\rho v_1) \\
e^{\mu(t-s)\D} \mathbb{P}(F_1(\rho_1,w_1,w_2)-F_1(\rho,v_1,v_2)) +e^{c_2(t-s)\D} \mathbb{Q}(F(\rho_1,w_1,w_2)-F(\rho,v_1,v_2))\\
e^{\mu(t-s)\D} \mathbb{P}(F_1(\rho_1,w_1,w_2)-F_1(\rho,v_1,v_2)) +e^{c_1(t-s)\D} \mathbb{Q}(F(\rho_1,w_1,w_2)-F(\rho,v_1,v_2)) \\
\end{array}
\right)\ ds.
\end{aligned}
\label{tech1a}
\end{equation}
with 
$$F_1(\rho,v_1,v_2)=\frac{1}{2}({\rm div}(\rho v_1\otimes v_2)+{\rm div}(\rho v_2\otimes v_1)).$$
Now we have:
\begin{equation}
\begin{aligned}
&{\rm div}(\rho v_1\otimes v_2)-{\rm div}(\rho_1 w_1\otimes w_2)=\\
&{\rm div}(\frac{1}{\rho}(\rho v_1-\rho_1 w_1)\otimes \rho v_2)+{\rm div}(\rho_1 w_1\otimes (\rho v_2-\rho_1 w_2)\frac{1}{\rho})+{\rm div}(\rho_1 w_1\otimes \rho_1 w_2\,\frac{(\rho_1-\rho)}{\rho\rho_1})
\end{aligned}
\label{tech1}
\end{equation}
From (\ref{estim1}) and (\ref{tech1}), we deduce that for $C,C_1>0$ large enough and for any $(\rho_1,\rho_1w_1,\rho_1w_2)$, $(\rho,\rho v_1,\rho v_2)$ in $ E_{R,M,T}$ we have for $C,C_1>0$ large enough:
\begin{equation}
\begin{aligned}
&\|\int^t_0 e^{\mu(t-s)\D} \mathbb{P}(F_1(\rho_1,\rho_1w_1,\rho_1w_2)-F_1(\rho,v_1,v_2)) ds\|_{{\cal E}_T}\\
&\leq C \|\frac{1}{\rho}\|_{L^\infty_T(L^\infty)}\big(\|\rho v_1-\rho_1 w_1\|_{{\cal E}_T}\|\rho v_2\|_{{\cal E}_T}+\|\rho v_2-\rho_1 w_2\|_{{\cal E}_T}\|\rho_1 w_1\|_{{\cal E}_T}\\
&+\|\rho_1 w_1\|_{{\cal E}_T} 
\|\rho_1 w_2\|_{{\cal E}_T} \|\frac{1}{\rho_1}\|_{L^\infty_T(L^\infty)}\|\rho-\rho_1\|_{L^\infty_T(L^\infty)}
+\|\rho v_2-\rho_1 w_2\|_{{\cal E}_T}\|\rho v_1\|_{{\cal E}_T}\\
&+\|\rho v_1-\rho_1 w_1\|_{{\cal E}_T}\|\rho_1 w_2\|_{{\cal E}_T}\big)\\
&\leq C_1M(R+M \beta R^2) \|(\rho,\rho v_1,\rho v_2)-(\rho_1,\rho_1w_1,\rho_1w_2)\|_{E_{R,M,T}}
\end{aligned}
\label{estimc1}
\end{equation}
We proceed similarly for the part $\|\int^t_0 e^{c_i(t-s)\D} \mathbb{Q}(F_1(\rho_1,\rho_1w_1,\rho_1w_2)-F_1(\rho,v_1,v_2)) ds\|_{{\cal E}_T}$ with $i=1,2$.
Similarly from (\ref{estim2}), we have for $c=c_1,c_2$ and $C>0$ large enough:
\begin{equation}
\begin{aligned}
&\|\int^t_0 e^{c(t-s)\D}\n (P(\rho)-P(\rho_1))ds\|_{{\cal E}_T}
\leq  CT\|P(\rho)-P(\rho_1)\|_{L^\infty_T(L^\infty(\R^N))}\\
&\leq C M_2 T\beta \|(\rho,\rho v_1,\rho v_2)-(\rho_1,\rho_1w_1,\rho_1w_2)\|_{E_{R,M,T}},
\end{aligned}
\label{estimc2}
\end{equation}
with $M_2=\sup_{x\in [0,2\|\rho_0\|_{L^\infty(\R^N)}]} |P'(x)|$. From (\ref{tech1a}), (\ref{tech1}), (\ref{estimc1}) and (\ref{estimc2})
we deduce that it exists $C>0$ sufficiently large such that:
\begin{equation}
\begin{aligned}
&\|\psi_2(\rho,\rho v_1,\rho v_2)-\psi_2(\rho_1,\rho_1 w_1,\rho_1 w_2),\psi_3(\rho,\rho v_1,\rho v_2)-\psi_3(\rho_1,\rho_1 w_1,\rho_1 w_2)\|_{{\cal E}_T}\\
&\leq C[M(R+M\beta R^2)+M_2 \beta T] \|(\rho,\rho v_1,\rho v_2)-(\rho_1,\rho_1w_1,\rho_1w_2)\|_{E_{R,M,T}}
\end{aligned}
\label{tech2}
\end{equation}
From (\ref{estim5}), it yields that for $C>0$ large enough:
\begin{equation}
\begin{aligned}
&\frac{1}{\beta}\|\psi_1(\rho,\rho v_1,\rho v_2)-\psi_1(\rho_1,\rho_1 w_1,\rho_1 w_2)\|_{L^\infty_T(L^\infty(\R^N))}\leq \frac{C}{\beta}\|(\rho v_1-\rho_1w_1)\|_{{\cal E}_T}
\end{aligned}
\label{tech3}
\end{equation}
Combining (\ref{tech2}) and (\ref{tech3}), we deduce that we have for $C>0$ large enough:
\begin{equation}
\begin{aligned}
&\|\psi(\rho,\rho v_1,\rho v_2)-\psi(\rho_1,\rho_1 w_1,\rho_1 w_2)\|_{{\cal E}_T}\\
&\leq C \big(M(R+M\beta R^2)+M_2  \beta T+\frac{1}{\beta} \big) \|(\rho,\rho v_1,\rho v_2)-(\rho_1,\rho_1w_1,\rho_1w_2)\|_{E_{R,M,T}}
\end{aligned}
\label{techfin}
\end{equation}
It suffices now to choose $\beta$, $T$ and $R$ such that:
\begin{equation}
 C \big(M(R+M\beta R^2)+M_2  \beta T+\frac{1}{\beta} \big) <1
 \label{cont}
 \end{equation}
 It proves in particular that the map $\psi$ is contractive and it concludes the proof of the theorem \ref{theo2}.
%
\\[2mm]
Let us prove now the Corollary \ref{cor1}, if  $\rho_0 v_1(0,\cdot)\in\overline{ ({\cal D}(\R^N))^N}^{(bmo^{-1}(\R^N))^N}$ then we can observe that:
\begin{equation}
\lim_{T\rightarrow 0}\|e^{t\D}\rho_0 v_1(0,\cdot)\|_{{\cal E}_T}=0.
\label{tech}
\end{equation}
Indeed assume that $w_0\in ({\cal D}(\R^N))^N$ then we have using the maximum principle for the heat equation and $C>0$ large enough:
$$
\begin{aligned}
\|e^{t\D}w_0 \|^2_{{\cal E}_T}&\leq \|w_0\|_{L^\infty} (\sqrt{T}+\sup_{0<t<T,x_0\in\R^N}\frac{1}{t^{\N}}\int^t_0\int_{B(x_0,\sqrt{t})} dx ds)\\
&\leq CT \|w_0\|_{L^\infty} .
\end{aligned}
$$
By density we can conclude. In particular it implies that for $T>0$ small enough we have:
\begin{equation}
\|e^{t\D}\rho_0 v_1(0,\cdot)\|_{{\cal E}_T}+\|e^{t\D}\rho_0 v_2(0,\cdot)\|_{{\cal E}_T}< \e_1.
\label{tech1v}
\end{equation}
Using the proof of the Theorem \ref{theo1}, we conclude that there exists a strong solution $(\rho,\rho u)$ on a finite time interval $[0,T]$. Indeed it we follow the previous proof, it suffices to fix $R$ and $T$ sufficiently small and verifying the previous estimate of the proof of the Theorem \ref{theo1}. Indeed by density it exist $w_1\in ({\cal D}(\R^N))^N$ such that:
$$\|e^{t\D} (w_1-\rho_0 v_1(0,\cdot))\|_{{\cal E}_T}\leq \|w_1-\rho_0 v_1(0,\cdot)\|_{bmo^{-1}_T}\leq \frac{R}{2}.$$
Now we choose $T_1\leq T$ sufficiently small such that $\|e^{t\D} w_1\|_{{\cal E}_{T_1}}\leq  \frac{R}{2}$, we get:
$$\|e^{t\D}(\rho_0 v_1(0,\cdot))\|_{{\cal E}_{T_1}}\leq R,$$
with $R$, $T_1$ satisfying the estimates of the proof of the Theorem \ref{theo1}.
 \section{Proof of the Theorem \ref{theo2}}
 In order to prove the Theorem \ref{theo2}, we are going to start by studying the linear system associated to (\ref{kfinal}):
 \begin{equation}
\begin{cases}
\begin{aligned}
&\p_t q-c\D q+{\rm div}m=F\\
&\p_t m-\mu\D m-\alpha\n{\rm div}m+\beta\n q=G,
\end{aligned}
\label{klinear}
\end{cases}
\end{equation}
with $c>0$, $\mu>0$, $\mu+\alpha>0$ and $\beta>0$. $(F,G)$ are external forces.
In the sequel we will note $W_{c,\mu,\alpha,\beta}$ the semigroup associated to the previous system and we have in particular from the Duhamel formula:
\begin{equation}
\begin{aligned}
&\left(\begin{array}{c}
q\\
m\\
\end{array}
\right)(t,\cdot)=W_{c,\mu,\alpha,\beta}(t)\left(\begin{array}{c}
q_0\\
m_0\\
\end{array}
\right)+\int_{0}^{t}W_{c,\mu,\alpha,\beta}(t-s)\left(\begin{array}{c}
 F \\
G\\
\end{array}
\right)(s)\ ds.
\end{aligned}
\end{equation}
This system has been studied by Bahouri et al (see \cite{BCD}) in the framework of the Besov space $B^{s}_{2,1}$ when $c=0$. We are going now to extend this study to the case of general Besov space of the form $B^{s}_{2,r}$. We set now $d=\Lambda^{-1}{\rm div}m$ and $\Omega=\Lambda^{-1}{\rm curl} m,$
with $\hat{\Lambda^{s_1} f }(\xi)=|\xi|^{s_1} \hat{f}(\xi)$ when $s_1\in\R$ and for $f$ a temperated distribution. We now study the following system:
\begin{equation}
\begin{cases}
\begin{aligned}
&\p_t q-c\D q+\Lambda d=F\\
&\p_t d-\nu\D d-\beta\Lambda q=\Lambda^{-1}{\rm div} G\\
&\p_t \Omega-\mu\D\Omega=\Lambda^{-1}{\rm curl} G.
\end{aligned}
\label{klinear1}
\end{cases}
\end{equation}
with $\nu=(\mu+\alpha)$. We refer to \cite{BCD}  for the definition of the Chemin-Lerner spaces $\widetilde{L}^\rho_T (B^{s}_{p,r})$ with $(\rho,p,r)\in[1,+\infty]^3$, $T>0$, $s\in\R$ and to \cite{arma} for the definition of the hybrid Besov spaces $\widetilde{B}^{s_1,s_2}_{p_1,p_2,r}$ with $(s_1,s_2)\in\R^2$,  $(p_1,p_2,r)\in[1,+\infty]^3$. 
\begin{proposition}
\label{flinear3}  
Let $T\in (0,+\infty]$. We assume that $(q_{0},u_{0})$ belongs to $\widetilde{B}_{2,\infty}^{s-1,s}\times (B_{2,\infty}^{s-1})^{N}$ with the source terms
$(F,G)$ in $\widetilde{L}^{1}_{T}(B^{s-1}_{2,\infty}\cap B^{s}_{2,\infty})\times
(\widetilde{L}^{1}_{T}(B^{s}_{2,\infty}))^{N}$.\\[1,5mm]
Let $(q,u)\in\big(\widetilde{L}^{\infty}_{T}(B^{s-1}_{2,\infty}\cap B^{s}_{2,\infty})\cap \widetilde{L}^{1}_{T}(B^{s+1}_{2,\infty}\cap B^{s+2}_{2,\infty})\big)\times\big(
 \widetilde{L}^{\infty}_{T}(B^{s-1}_{2,\infty})\cap \widetilde{L}^{1}_{T}(B^{s+1}_{2,\infty})\big)^{N}$ be a solution of the
system $(\ref{klinear})$, then there exists a universal constant $C$ such that for any $T>0$ we have:
\begin{equation}
\begin{aligned}
&\|(q,\n q,u)\|_{\ \widetilde{L}^{\infty}_{T}(B^{s-1}_{2,\infty})\cap \widetilde{L}^{1}_{T}(B^{s+1}_{2,\infty})}\leq C(\|(q_0,\n q_{0},u_{0})\|_{B^{s-1}_{2,\infty}}+\|(F,\n
F,G)\|_{ \widetilde{L}^{1}_{T}(B^{s-1}_{2,\infty})}).
\end{aligned}
\label{lineairestim}
\end{equation}
\end{proposition}
{\bf Proof:}
Let $(q,m)$ be a solution of $(\ref{klinear})$, we are going to separate the case of the low and high frequencies,
which have a
different behavior concerning the control of the derivative index for the Besov spaces. Our goal consists now in studying the system (\ref{klinear1}) and in particular to estimate $(q,d)$ and $\Omega$. We observe that $\Omega$ verifies simply an heat equation and classical estimates on the heat equation in Besov spaces give (see \cite{BCD}):
\begin{equation}
\|\Omega\|_{\widetilde{L}^{\infty}_{T}(B^{s-1}_{2,\infty})\cap \widetilde{L}^{1}_{T}(B^{s+1}_{2,\infty})}\leq C(\|\Omega_{0}\|_{B^{s-1}_{2,\infty}}+\|G\|_{ \widetilde{L}^{1}_{T}(B^{s-1}_{2,\infty})}).
\label{curl}
\end{equation}
Let us study now the unknowns $(q,d)$.
\subsubsection*{Case of low frequencies}
We assume here that $l\in\mathbb{Z}$ with $l\leq l_0$ (we will determine later $l_0\in\mathbb{Z}$).
Applying operator $\D_{l}$ (see \cite{BCD} for the definition of $\D_l$) to the system (\ref{klinear1}) and denoting $g_l=\D_l g$, we obtain
the following system:
\begin{equation}
\begin{cases}
\begin{aligned}
&\frac{d}{dt}q_{l}+c\Lambda^2 q_l+\Lambda
d_{l}=F_{l},\\[2mm]
&\frac{d}{dt}d_{l}+\nu\Lambda^2
d_{l}-\beta \Lambda q_{l}=G_{l}.
\end{aligned}
\end{cases}
\label{5LPH}
\end{equation}
We set:
\begin{equation}
f_{l}^{2}=\beta \|q_{l}\|_{L^{2}}^{2} +\|d_{l}
\|_{L^{2}}^{2}.
\label{5P0}
\end{equation}
Taking the $L^{2}$ scalar product of the first equation of
(\ref{5LPH}) with $q_{l}$ and of the second equation with
$d_{l}$, we get the following two identities:
\begin{equation}
\begin{cases}
\begin{aligned}
&\frac{1}{2}\frac{d}{dt}\|q_{l}\|_{L^{2}}^{2}+(\Lambda d_{l},q_{l})+c\|\Lambda q_l\|_{L^2}^2=(F_{l},q_{l}),\\[2mm]
&\frac{1}{2}\frac{d}{dt}\|d_{l}\|_{L^{2}}^{2}+\nu\|\Lambda
d_{l}\|_{L^{2}}^{2} -\beta(\Lambda q_{l},d_{l})
=(G_{l},d_{l}).
\end{aligned}
\end{cases}
\label{5P1}
\end{equation}
We deduce that:
\begin{equation}
\begin{aligned}
&\frac{1}{2}\frac{d}{dt} f_l^2+(\beta c\|\Lambda q_l\|_{L^2}+\nu\|\Lambda
d_{l}\|_{L^{2}}^{2})
\leq \|G_{l}\|_{L^2}\|d_{l}\|_{L^2}+\beta\|F_l\|_{L^2}\|q_{l}\|_{L^2}.
\end{aligned}
\label{5P2}
\end{equation}
From the definition of $f_l^2$ we deduce that it exists $\alpha_1,C>0$ independent on $l_0$ such that:
\begin{equation}
\frac{1}{2}\frac{d}{dt}f_{l}^{2}+\alpha_1 2^{2l}f_{l}^{2}\leq Cf_l(\|G_{l}\|_{L^2}+\|F_l\|_{L^2}).
\label{5P5}
\end{equation}
We have in particular for $l\leq l_0$:
\begin{equation}
\begin{aligned}
&\frac{1}{2}\frac{d}{dt}f_{l}+\alpha_1 2^{2l}f_{l}\leq C(\|G_{l}\|_{L^2}+\|F_l\|_{L^2}).
\end{aligned}
\label{5P6}
\end{equation}
\subsubsection*{Case of high frequencies}
We consider now the case where $l\geq l_{0}+1$ and we define now $f_l$ as follows:
$$f_{l}^{2}=\|\Lambda q_{l}\|_{L^{2}}^{2}+A\|d_{l}\|_{L^{2}}^{2}-
\frac{2}{c+\nu}(\Lambda q_{l},d_{l}),$$
with $A>0$ to be determinated. We apply the
operator $\Lambda\D_l$ to the first equation of (\ref{klinear1}), multiply
by $\Lambda q_{l}$ and integrate over $\R^{N}$, so we obtain:
\begin{equation}
\frac{1}{2}\frac{d}{dt}\|\Lambda
q_{l}\|_{L^{2}}^{2}+c\|\Lambda^2 q_l\|_{L^2}^2+(\Lambda^{2}d_{l},\Lambda q_{l})=(\Lambda F_l,\Lambda q_l).
 \label{5H7}
\end{equation}
Moreover we have in a similar way:
\begin{equation}
\begin{cases}
\begin{aligned}
&\frac{1}{2}\frac{d}{dt}\|d_{l}\|_{L^{2}}^{2}+\bar{\nu}\|\Lambda
d_{l}\|_{L^{2}}^{2} -\beta(\Lambda q_{l},d_{l})
=(G_{l},d_{l}).
\\
&\frac{d}{dt}(\Lambda q_{l},d_{l})+\|\Lambda d_{l}\|_{L^{2}}^{2}-
\beta \|\Lambda q_{l}\|_{L^{2}}^{2}+(c+\nu)(\Lambda^{2}d_{l}, \Lambda q_{l})
=(\Lambda F_l,d_l)+(G_l,\Lambda q_l).\\
\end{aligned}
\end{cases}
\label{5H8}
\end{equation}
By linear combination of (\ref{5H7})-(\ref{5H8}) we have:
\begin{equation}
\begin{aligned}
&\frac{1}{2}\frac{d}{dt}f_{l}^{2}+c\|\Lambda^2
q_{l}\|_{L^{2}}^{2}+\big(A\nu-\frac{1}{c+\nu}\big)\|\Lambda
d_{l}\|_{L^{2}}^{2}+\frac{\beta}{c+\nu}\|\Lambda q_l\|_{L^2}^2-A\beta(\Lambda
q_{l},d_{l})\\
&=-\frac{1}{c+\nu}(\Lambda F_l,d_l)-\frac{1}{c+\nu}(G_l,\Lambda q_l)+(\Lambda F_l,\Lambda q_l)+A(G_l,d_l).
\end{aligned}
\label{5P10}
\end{equation}
We have now in using Young inequalities for all $a>0$:
$$
\begin{aligned}
&|(d_{l},\Lambda q_{l})|\leq\frac{a}{2}\|\Lambda q_{l}\|_{L^{2}}^{2}+\frac{1}{2a}\|d_{l}\|_{L^{2}}^{2},\\
\end{aligned}
$$
\begin{equation}
\begin{aligned}
&\frac{1}{2}\frac{d}{dt}f_{l}^{2}+c\|\Lambda^2
q_{l}\|_{L^{2}}^{2}+\big(A\nu-\frac{1}{c+\nu}\big)\|\Lambda
d_{l}\|_{L^{2}}^{2}+(\frac{\beta}{c+\nu}-\frac{Aa\beta}{2})\|\Lambda q_l\|_{L^2}^2-\frac{A\beta}{2a}\|d_l\|_{L^2}^2\\
&\leq -\frac{1}{c+\nu}(\Lambda F_l,d_l)-\frac{1}{c+\nu}(G_l,\Lambda q_l)+(\Lambda F_l,\Lambda q_l)+A(G_l,d_l).
\end{aligned}
\label{5P12}
\end{equation}
We have now since $l\geq l_0+1$
 for $c_1>0$:
\begin{equation}
\begin{aligned}
&\frac{1}{2}\frac{d}{dt}f_{l}^{2}+c\|\Lambda^2
q_{l}\|_{L^{2}}^{2}+\big(A\nu-\frac{1}{c+\nu}-2^{-2l _0}\frac{A c_1 \beta}{2a}\big)\|\Lambda
d_{l}\|_{L^{2}}^{2}+(\frac{\beta}{c+\nu}-\frac{Aa\beta}{2})\|\Lambda q_l\|_{L^2}^2\\
&\leq -\frac{1}{c+\nu}(\Lambda F_l,d_l)-\frac{1}{c+\nu}(G_l,\Lambda q_l)+(\Lambda F_l,\Lambda q_l)+A(G_l,d_l).
\end{aligned}
\label{5P13}
\end{equation}
We choose now $a$,
$A$ and $l_0$ as follows:
\begin{equation}
\begin{aligned}
&A \nu=\frac{2 M}{c+\nu},\;\frac{\beta}{c+\nu}-\frac{Aa\beta}{2} \geq 0,\;A\nu-\frac{1}{c+\nu}-2^{-2l _0}\frac{A\beta}{2a}\geq \frac{1}{2(c+\nu)}.
\end{aligned}
\end{equation}
with $M>1$ sufficiently large to determine later. Now from the definition of $f_l^2$ and using Young inequality, we have:
\begin{equation}
\frac{1}{2} \|\Lambda q_l\|_{L^2}^2+ (\frac{2M}{\nu(c+\nu)}-\frac{2}{(c+\nu)^2})\|d_l\|_{L^2}^2\leq f_l^2\leq \frac{3}{2} \|\Lambda q_l\|_{L^2}^2+ (\frac{2M}{\nu(c+\nu)}+\frac{2}{(c+\nu)^2})\|d_l\|_{L^2}^2
\label{high}
\end{equation}
We choose now $M$ sufficiently large such that $\frac{2M}{\nu(c+\nu)}-\frac{2}{(c+\nu)^2}\geq\frac{1}{2}$.
Using (\ref{5P13}) and (\ref{high}) we deduce that there exists  constants $\alpha_2>0$ and $C_1>0$ such that for $l\geq l_{0}+1$ we
have:
\begin{equation}
\frac{1}{2}\frac{d}{dt}f_{l}^{2}+\alpha_2 2^{2l} f_{l}^{2}\leq C_1 f_l(2^l\|F_l\|_{L^2}+\|G_l\|_{L^2}).
\label{5P14}
\end{equation}
It yields that:
\begin{equation}
\frac{1}{2}\frac{d}{dt}f_{l}+\alpha_2 2^{2l} f_{l}\leq C_1 (2^l\|F_l\|_{L^2}+\|G_l\|_{L^2}).
\label{5P15}
\end{equation}
\subsubsection*{Final estimates}
Integrating over $[0,t]$ the estimates (\ref{5P6}) and (\ref{5P15}) and multiplying by $2^{l(s-1)}$, we have for $C_2>0$ large enough and any $l\in\mathbb{Z}$:
\begin{equation}
\begin{aligned}
&\frac{1}{2}2^{l(s-1)}f_{l}(t)+\min(\alpha_1,\alpha_2) 2^{2ls}\int^t_0 f_{l}(s) ds \\
&\hspace{3cm}\leq \frac{1}{2}f_l(0)+C_2(\|G_{l}\|_{\widetilde{L}^1_t(B^{s-1}_{2,\infty})}+\|F_l\|_{\widetilde{L}^1_t(B^{s-1}_{2,\infty}\cap B^{s-1}_{2,\infty})}).
\end{aligned}
\label{estimofin}
\end{equation}
From the definition of $f_l$ we deduce the estimate (\ref{lineairestim}) using in particular (\ref{high}). It concludes the proof of the proposition.
\hfill {$\Box$}\\
\\
Let us study again the system (\ref{klinear}) and applying ${\rm div}$ to the momentum equation, we have:
 \begin{equation}
\begin{cases}
\begin{aligned}
&\p_t q-c\D q+{\rm div}m=0\\
&\p_t {\rm div}m-\nu \D {\rm div}m+\beta\D  q=0,
\end{aligned}
\label{klinearp}
\end{cases}
\end{equation}
We assume now that $\mu, c,\nu,\beta>0$ and $\nu\ne c$. The only case where $\nu=c$ will be in the sequel the case $\kappa^2=\mu^2$. We will study this case later.\\
When we apply the Fourier transform ${\cal F}$, we have then:
\begin{equation}
\begin{aligned}
{\cal F}&\left(\begin{array}{c}
q\\
{\rm div}m\\
\end{array}
\right)(t,\xi)=e^{-t A(\xi)}{\cal F}\left(\begin{array}{c}
q_0\\
{\rm div}m_0\\
\end{array}
\right)(\xi)
\end{aligned}
\end{equation}
with:
$
A(\xi)=\quad
\begin{pmatrix} 
c|\xi|^2 & 1 \\
-\beta|\xi|^2 & \nu|\xi|^2 
\end{pmatrix}
\quad$.
The characteristic polynom is:
$$P_{A(\xi)}(\lambda)=\lambda^2-\lambda|\xi|^2(c+\nu)+c\nu|\xi|^4+\beta|\xi|^2.$$
The eigenvalues are:
\begin{itemize}
\item  if $|\xi|^2 \geq \frac{4\beta}{(\nu-c)^2}$:
$$
\begin{aligned}
&\lambda_{1h}(\xi)=\frac{1}{2}|\xi|^2\big((\nu+c)|+|\nu-c| \sqrt{1-\frac{4\beta}{(\nu-c)^2|\xi|^2}}\big),\\&\lambda_{2h}(\xi)=\frac{1}{2}|\xi|^2\big((\nu+c)|-|\nu-c| \sqrt{1-\frac{4\beta}{(\nu-c)^2|\xi|^2}}\big).
\end{aligned}
$$
\item if $|\xi|^2< \frac{4\beta}{(\nu-c)^2}$:
$$
\begin{aligned}
&\lambda_{1l}(\xi)=\frac{1}{2}|\xi|^2\big((\nu+c)|+i|\nu-c| \sqrt{\frac{4\beta}{(\nu-c)^2|\xi|^2}-1}\big),\\&\lambda_{2l}(\xi)=\frac{1}{2}|\xi|^2\big((\nu+c)|-i|\nu-c| \sqrt{\frac{4\beta}{(\nu-c)^2|\xi|^2}-1}\big).
\end{aligned}
$$
\end{itemize}
\subsection*{High Frequencies}
When $|\xi|^2 \geq \frac{4\beta}{(\nu-c)^2}$ we have: 
\begin{equation}
\begin{aligned}
&{\cal F}&\left(\begin{array}{c}
q\\
{\rm div}m\\
\end{array}
\right)(t,\xi)=P(\xi) 
\begin{pmatrix} 
e^{-t\lambda_{1h}(\xi)} & 0 \\
0&e^{-t\lambda_{2h}(\xi)}
\end{pmatrix} (P(\xi))^{-1}
 {\cal F}\left(\begin{array}{c}
q_0\\
{\rm div}m_0\\
\end{array}
\right)(\xi),
\end{aligned}
\end{equation}
with:
$$
\begin{aligned}
&P(\xi)=\begin{pmatrix} 
\frac{-2}{|\xi|^2\big(c-\nu-|c-\nu|\alpha(\xi)\big)}  & \frac{-2}{|\xi|^2\big(c-\nu+|c-\nu|\alpha(\xi)\big)} \\
1 & 1
\end{pmatrix}
\\
&
(P(\xi))^{-1}=-\frac{\beta}{|c-\nu|\alpha(\xi)}
\begin{pmatrix} 
 1& \frac{2}{|\xi|^2\big(c-\nu+|c-\nu|\alpha(\xi)\big)} \\
-1 & \frac{-2}{|\xi|^2\big(c-\nu-|c-\nu|\alpha(\xi)\big)} 
\end{pmatrix}
.
\end{aligned}
$$
We have noted $\alpha(\xi)=\sqrt{1-\frac{4\beta}{(\nu-c)^2|\xi|^2}}$.
Finally we get:
\begin{equation}
\begin{aligned}
&{\cal F}\left(\begin{array}{c}
q\\
m\\
\end{array}
\right)(t,\xi)=-\frac{\beta}{|c-\nu|\alpha(\xi)} \\
&\begin{pmatrix} 
 \frac{-2 e^{-t\lambda_{1h}(\xi)}}{|\xi|^2\big(c-\nu-|c-\nu|\alpha(\xi)\big)}  
 + \frac{2 e^{-t\lambda_{2h}(\xi)}}{|\xi|^2\big(c-\nu+|c-\nu|\alpha(\xi)\big)}  
  &\frac{1}{\beta|\xi|^2}(e^{-t\lambda_{2h}(\xi)}-e^{-t\lambda_{1h}(\xi)})\\
e^{-t\lambda_{1h}(\xi)}-e^{-t\lambda_{2h}(\xi)} & \frac{2 e^{-t\lambda_{1h}(\xi)}}{|\xi|^2\big(c-\nu+|c-\nu|\alpha(\xi)\big)}  
 - \frac{2 e^{-t\lambda_{2h}(\xi)}}{|\xi|^2\big(c-\nu-|c-\nu|\alpha(\xi)\big)}  \\
\end{pmatrix}\\
&\times {\cal F}\left(\begin{array}{c}
q_0\\
{\rm div}m_0\\
\end{array}
\right)(\xi).
\end{aligned}
\end{equation}
It yields then:
\begin{equation}
\begin{aligned}
&{\cal F}
q (t,\xi) =(e^{-t\lambda_{1h}(\xi)}-e^{-t\lambda_{2h}(\xi)})\big(  \frac{\mbox{sgn}(c-\nu)}{2\alpha(\xi)}  {\cal F}q_0(\xi)+\frac{1}{|\xi|^2|c-\nu| \alpha(\xi)}{\cal F}{\rm div}m_0(\xi)\big)\\
&+\frac{1}{2}(e^{-t\lambda_{2h}(\xi)}+e^{-t\lambda_{1h}(\xi)})  {\cal F}q_0(\xi)\\
&{\cal F}
{\rm div}m(t,\xi) =(e^{-t\lambda_{2h}(\xi)}-e^{-t\lambda_{1h}(\xi)})\big( \frac{\beta}{|c-\nu|\alpha(\xi)}{\cal F}q_0(\xi)+\frac{\mbox{sgn}(c-\nu)}{2\alpha(\xi)}{\cal F}{\rm div}m_0(\xi)\big)\\
&+\frac{1}{2}(e^{-t\lambda_{2h}(\xi)}+e^{-t\lambda_{1h}(\xi)}){\cal F}{\rm div}m_0(\xi).
\end{aligned}
\label{finh1}
\end{equation}
\subsection*{Low Frequencies}
When $|\xi|^2 < \frac{4\beta}{(\nu-c)^2}$ we have:
\begin{equation}
\begin{aligned}
&{\cal F}&\left(\begin{array}{c}
q\\
{\rm div}m\\
\end{array}
\right)(t,\xi)=P_1(\xi) 
\begin{pmatrix} 
e^{-t\lambda_{1l}(\xi)} & 0 \\
0&e^{-t\lambda_{2l}(\xi)}
\end{pmatrix} (P_1(\xi))^{-1}
{\cal F}\left(\begin{array}{c}
q_0\\
{\rm div} m_0\\
\end{array}
\right)(\xi).
\end{aligned}
\end{equation}
In addition if we set $\alpha_1(\xi)=\sqrt{\frac{4\beta}{(\nu-c)^2|\xi|^2}-1}$
, we have:
$$
\begin{aligned}
&P_1(\xi)=\begin{pmatrix} 
\frac{-2}{|\xi|^2\big(c-\nu-i|c-\nu|\alpha_1(\xi)\big)}  & \frac{-2}{|\xi|^2\big(c-\nu+i|c-\nu|\alpha_1(\xi)\big)} \\
1 & 1
\end{pmatrix}
\\
&(P_1(\xi))^{-1}=\frac{\beta i}{|c-\nu|\alpha_1(\xi)}
\begin{pmatrix} 
 1& \frac{2}{|\xi|^2\big(c-\nu+i|c-\nu|\alpha_1(\xi)\big)} \\
-1 & \frac{-2}{|\xi|^2\big(c-\nu-i|c-\nu|\alpha_1(\xi)\big)} 
\end{pmatrix}.
\end{aligned}
$$
Finally we obtain:
\begin{equation}
\begin{aligned}
&{\cal F}\left(\begin{array}{c}
q\\
m\\
\end{array}
\right)(t,\xi)=\frac{\beta i}{|c-\nu|\alpha_1 (\xi)} \\
&\begin{pmatrix} 
 \frac{-2 e^{-t\lambda_{1l}(\xi)}}{|\xi|^2\big(c-\nu-i |c-\nu|\alpha_1(\xi)\big)}  
 + \frac{2 e^{-t\lambda_{2l}(\xi)}}{|\xi|^2\big(c-\nu+i |c-\nu|\alpha_1(\xi)\big)}  
  &\frac{1}{\beta|\xi|^2}(e^{-t\lambda_{2l}(\xi)}-e^{-t\lambda_{1l}(\xi)})\\
e^{-t\lambda_{1l}(\xi)}-e^{-t\lambda_{2l}(\xi)} & \frac{2 e^{-t\lambda_{1l}(\xi)}}{|\xi|^2\big(c-\nu+i |c-\nu|\alpha_1(\xi)\big)}  
 - \frac{2 e^{-t\lambda_{2l}(\xi)}}{|\xi|^2\big(c-\nu-i |c-\nu|\alpha_1(\xi)\big)}  \\
\end{pmatrix}\\
&\times {\cal F}\left(\begin{array}{c}
q_0\\
{\rm div} m_0\\
\end{array}
\right)(\xi).
\end{aligned}
\end{equation}
We get then:
\begin{equation}
\begin{aligned}
&{\cal F}
q (t,\xi) =i (e^{-t\lambda_{2l}(\xi)}-e^{-t\lambda_{1l}(\xi)})\big(  \frac{\mbox{sgn}(c-\nu)}{2\alpha_1(\xi)}  {\cal F}q_0(\xi)+\frac{1}{|\xi|^2|c-\nu| \alpha_1(\xi)}{\cal F}{\rm div}m_0(\xi)\big)\\
&+\frac{1}{2}(e^{-t\lambda_{2l}(\xi)}+e^{-t\lambda_{1l}(\xi)}) {\cal F}q_0(\xi)\\
&{\cal F}
{\rm div}m(t,\xi) =i(e^{-t\lambda_{1l}(\xi)}-e^{-t\lambda_{2l}(\xi)})\big( \frac{\beta }{|c-\nu|\alpha_1(\xi)}{\cal F}q_0(\xi)+\frac{\mbox{sgn}(c-\nu)}{2\alpha_{1}(\xi)}{\cal F}{\rm div}m_0(\xi)\big)\\
&+\frac{1}{2}(e^{-t\lambda_{2l}(\xi)}+e^{-t\lambda_{1l}(\xi)}){\cal F}{\rm div}m_0(\xi).
\end{aligned}
\label{finh2}
\end{equation}
In the two next propositions, in order to simplify the notations we will denote by $W(t)$ the semigroup $W_{c,\mu,\alpha,\beta}(t)$ with $t>0$.
\begin{proposition}
\label{propocru}
Let $\phi$ be a smooth function supported in the shell ${\cal C}(0,R_1,R_2)$ with $0<R_1<R_2$. There exist two positive constants $\kappa$ and $C$ depending only on $\phi$ and such that for all $t\geq 0$ and $\lambda>0$, we have:
\begin{equation}
\|\phi(\lambda^{-1}D)W(t)\left(\begin{array}{c}
q_0\\
 m_0\\
\end{array}
\right)\|_{L^2}\leq C e^{-\kappa t \lambda^2}\|\phi(\lambda^{-1}D)\left(\begin{array}{c}
q_0\\
 m_0\\
\end{array}
\right)\|_{L^2}
\label{crucial}
\end{equation}
\end{proposition}
{\bf Proof:} From (\ref{finh1}), we deduce that for $|\xi|^2\geq  \frac{4\beta}{(\nu-c)^2}$ (with $R_k$ the Riez transform such that ${\cal F}R_kg (\xi)=\frac{i\xi_k}{|\xi|} {\cal F} g(\xi)$ for $g$ a temperated distribution):
\begin{equation}
\begin{aligned}
&\big[{\cal F} \phi(\lambda^{-1}D)W(t)\left(\begin{array}{c}
q_0\\
 m_0\\
\end{array}
\right)\big]_1 (t,\xi) \\
&=\phi (\frac{|\xi|}{\lambda})\biggl((e^{-t\lambda_{1h} (\xi)}-e^{-t\lambda_{2h}(\xi)})\big(  \frac{\mbox{sgn}(c-\nu)}{2\alpha(\xi)}  {\cal F}q_0(\xi)+\frac{1}{|\xi||c-\nu| \alpha(\xi)}\sum_{j=1}^N {\cal F}R_j (m_0)_j (\xi)\big)\\
&+\frac{1}{2}(e^{-t\lambda_{2h}(\xi)}+e^{-t\lambda_{1h}(\xi)})  {\cal F}q_0(\xi)\biggl)
\\[2mm]
&\big[{\cal F} \phi(\lambda^{-1}D)W(t)\left(\begin{array}{c}
q_0\\
 m_0\\
\end{array}
\right)\big]_{1+i} (t,\xi) \\
&=\phi (\frac{\xi}{\lambda})\biggl( (e^{-t\lambda_{2h}(\xi)}-e^{-t\lambda_{1h}(\xi)})\big( -\frac{\beta}{|\xi| |c-\nu|\alpha(\xi)}{\cal F}R_i q_0(\xi)+\frac{\mbox{sgn}(c-\nu)}{2\alpha(\xi)}{\cal F}(\mathbb{Q}m_0)_i (\xi)\big)\\
&+\frac{1}{2}(e^{-t\lambda_{2h}(\xi)}+e^{-t\lambda_{1h}(\xi)}){\cal F}(\mathbb{Q}m_0)_i(\xi)+e^{-\mu t|\xi|^2}{\cal F}(\mathbb{P}m_0)_i (\xi)\biggl),
\end{aligned}
\label{finh1a}
\end{equation}
with $i\in[1,N]$.\\ We have used the fact that $\mathbb{P}m$ satisfies simply an heat equation.
Similarly from (\ref{finh2}), we have for $|\xi|^2< \frac{4\beta}{(\nu-c)^2}$:
\begin{equation}
\begin{aligned}
&\big[{\cal F} \phi(\lambda^{-1}D)W(t)\left(\begin{array}{c}
q_0\\
 m_0\\
\end{array}
\right)\big]_1 (t,\xi)  =\phi (\frac{|\xi|}{\lambda})\biggl(i (e^{-t\lambda_{2l}(\xi)}-e^{-t\lambda_{1l}(\xi)})\big(  \frac{\mbox{sgn}(c-\nu)}{2\alpha_1(\xi)}  {\cal F}q_0(\xi)\\
&+\frac{1}{|\xi||c-\nu| \alpha_1(\xi)}\sum_{j=1}^N {\cal F}R_j (m_0)_j (\xi) \big)+\frac{1}{2}(e^{-t\lambda_{2l}(\xi)}+e^{-t\lambda_{1l}(\xi)})  {\cal F}q_0(\xi)\biggl)
\\[1,5mm]
&\big[{\cal F} \phi(\lambda^{-1}D)W(t)\left(\begin{array}{c}
q_0\\
 m_0\\
\end{array}
\right)\big]_{1+i} (t,\xi) =\phi (\frac{\xi}{\lambda})\biggl(i(e^{-t\lambda_{1l}(\xi)}-e^{-t\lambda_{2l}(\xi)})\big( \frac{-\beta }{|c-\nu||\xi|\alpha_1(\xi)}{\cal F}R_iq_0(\xi)\\
&+\frac{\mbox{sgn}(c-\nu)}{2\alpha_1(\xi)}{\cal F}(\mathbb{Q}m_0)_i (\xi)\big)+\frac{1}{2}(e^{-t\lambda_{2l}(\xi)}+e^{-t\lambda_{1l}(\xi)})\big){\cal F}(\mathbb{Q}m_0)_i  (\xi)+e^{-\mu t|\xi|^2}{\cal F}(\mathbb{P}m_0)_i (\xi)\biggl).
\end{aligned}
\label{finh2a}
\end{equation}
Applying Plancherel Theorem, we observe easily that it exists $\e>0$ small enough and $C_\e,\kappa_\e>0$ such that:
\begin{equation}
\begin{aligned}
& \| 1_{\R^N / {\cal C}(0,\sqrt{\frac{4\beta}{|\nu-c|^2}-\e},\sqrt{\frac{4\beta}{|\nu-c|^2}+\e})  }{\cal F} \phi(\lambda^{-1}D)W(t)\left(\begin{array}{c}
q_0\\
 m_0\\
\end{array}
\right)\|_{L^2}\\
&\hspace{7cm}\leq C_\e e^{-\kappa _\e\lambda^2 t}\|\phi(\lambda^{-1}D)\left(\begin{array}{c}
q_0\\
 m_0\\
\end{array}
\right)\|_{L^2}
\end{aligned}
\label{estim3u}
\end{equation}
Indeed we use the fact that when $\xi\in \R^N / {\cal C}(0,\sqrt{\frac{4\beta}{|\nu-c|^2}-\e},\sqrt{\frac{4\beta}{|\nu-c|^2}+\e})$ we have:
$$
\begin{aligned}
&\alpha_1(\xi)\geq |\nu-c|\sqrt{\frac{\e}{4\beta(1-\frac{\e(\nu-c)^2}{4\beta})}},\;\alpha(\xi)\geq  |\nu-c|\sqrt{\frac{\e}{4\beta(1+\frac{\e(\nu-c)^2}{4\beta})}},\; \alpha_1(\xi)|\xi|\geq \sqrt{\e}.
\end{aligned}
$$
The only difficulty is the behavior of the solution in the region ${\cal C}(0,\sqrt{\frac{4\beta}{|\nu-c|^2}-\e},\sqrt{\frac{4\beta}{|\nu-c|^2}+\e}) $, in particular when ${\cal C}(0,R_1\lambda,R_2\lambda)\cap {\cal C}(0,\sqrt{\frac{4\beta}{|\nu-c|^2}-\e},\sqrt{\frac{4\beta}{|\nu-c|^2}+\e})\ne \emptyset$. We have in particular when $|\xi|^2 \in (\frac{4 \beta}{|\nu-c|^2},\frac{4\beta }{|\nu-c|^2}+\e)$ that $\alpha(\xi)\in (0,\sqrt{\e})$:
$$
\begin{aligned}
&\big|\frac{e^{-t\lambda_1(\xi)}-e^{-t\lambda_2(\xi)}}{\alpha(\xi)}\big|=e^{-\frac{t|\xi|^2(\nu+c)}{2}}\big|\frac{e^{-\frac{t|\xi|^2|\nu-c|\alpha(\xi)}{2}}-e^{\frac{t|\xi|^2|\nu-c|\alpha(\xi)}{2}}}{\alpha(\xi)}\big|
\end{aligned}
$$
When $t\alpha(\xi)\leq 1$, it exists $C>0$ sufficiently large such that for any $x\in [0,|\nu-c|(\frac{4\beta }{|\nu-c|^2}+\e)]$ we have:
$$e^{x}-e^{-x}\leq C x.$$
It implies that we have when $t\alpha(\xi)\leq 1$, it exists $C_1,C_2>0$ large enough such that:
\begin{equation}
\begin{aligned}
&\big|\frac{e^{-t\lambda_1(\xi)}-e^{-t\lambda_2(\xi)}}{\alpha(\xi)}\big|\leq C_1 t|\xi|^2 e^{-\frac{t|\xi|^2(\nu+c)}{2}}\leq C_2 e^{-\frac{t|\xi|^2(\nu+c)}{4}}.
\end{aligned}
\label{estim1u}
\end{equation}
When $t\alpha(\xi)\geq 1$ (and $\frac{1}{\alpha(\xi)}\leq t$), we deduce that for $C_3>0$ large enough and $\e\leq\frac{1}{2}$:
\begin{equation}
\begin{aligned}
\big|\frac{e^{-t\lambda_1(\xi)}-e^{-t\lambda_2(\xi)}}{\alpha(\xi)}\big|&\leq \frac{|\nu-c|^2}{4\beta}t |\xi|^2e^{-\frac{t|\xi|^2(\nu+c)}{2}}\big|e^{-\frac{t|\xi|^2|\nu-c|\alpha(\xi)}{2}}-e^{\frac{t|\xi|^2|\nu-c|\alpha(\xi)}{2}}\big|\\
&\leq \frac{|\nu-c|^2}{4\beta}t |\xi|^2e^{-\frac{t|\xi|^2(\nu+c)}{2}}+\frac{|\nu-c|^2}{4\beta}t |\xi|^2e^{-\frac{t|\xi|^2(\nu+c-\sqrt{\e}|\nu-c|)}{2}}\\
&\leq C_3 e^{-\frac{t|\xi|^2(\nu+c)}{4}}.
\end{aligned}
\label{estim2u}
\end{equation}
From (\ref{finh1a}), (\ref{estim3u}), (\ref{estim1u}) and (\ref{estim2u}), we deduce that
 it exists $\e>0$ small enough and $C_\e,\kappa_\e>0$ such that:
\begin{equation}
\begin{aligned}
& \| 1_{\R^N / {\cal C}(0,\frac{2\sqrt{\beta}}{|\nu-c|},\sqrt{\frac{4\beta}{|\nu-c|^2}+\e)}  }{\cal F} \phi(\lambda^{-1}D)W(t)\left(\begin{array}{c}
q_0\\
 m_0\\
\end{array}
\right)\|_{L^2}\leq C_\e e^{-\kappa _\e\lambda^2 t}\|\phi(\lambda^{-1}D)\left(\begin{array}{c}
q_0\\
 m_0\\
\end{array}
\right)\|_{L^2}
\end{aligned}
\label{estim4u}
\end{equation}
Let us deal now with the case $|\xi|^2\in (  \frac{4\beta}{|\nu-c|^2},\frac{4\beta}{|\nu-c|^2}-\e)$ whoch corresponds to $\alpha_1(\xi)\in(0,\sqrt{\e})$, we have then:
$$
\begin{aligned}
&\big|\frac{e^{-t\lambda_1(\xi)}-e^{-t\lambda_2(\xi)}}{\alpha_1(\xi)}\big|=e^{-\frac{t|\xi|^2(\nu+c)}{2}}\big|\frac{e^{-\frac{it|\xi|^2|\nu-c|\alpha_1(\xi)}{2}}-e^{\frac{it|\xi|^2|\nu-c|\alpha(\xi)}{2}}}{\alpha_1(\xi)}\big|\\
&=2 e^{-\frac{t|\xi|^2(\nu+c)}{2}}\frac{|\sin\big(\frac{t|\xi|^2|\nu-c|\alpha_1(\xi)}{2}\big)|}{|\alpha_1(\xi)|}.
\end{aligned}
$$
When $t\alpha_1(\xi)\leq 1$, it exists $C_1,C_2>0$ large enough such that:
\begin{equation}
\begin{aligned}
&\big|\frac{e^{-t\lambda_1(\xi)}-e^{-t\lambda_2(\xi)}}{\alpha_1(\xi)}\big|\leq C_1 t|\xi|^2 e^{-\frac{t|\xi|^2(\nu+c)}{2}}\leq C_2 e^{-\frac{t|\xi|^2(\nu+c)}{4}}.
\end{aligned}
\label{estim1uc}
\end{equation}
When $t\alpha_1(\xi)\geq 1$, we deduce that for $C_3>0$ large enough:
\begin{equation}
\begin{aligned}
\big|\frac{e^{-t\lambda_1(\xi)}-e^{-t\lambda_2(\xi)}}{\alpha_1(\xi)}\big|&\leq\frac{2}{ \frac{4\beta}{|\nu-c|^2}-\e} t |\xi|^2e^{-\frac{t|\xi|^2(\nu+c)}{2}}\\
&\leq C_3 e^{-\frac{t|\xi|^2(\nu+c)}{4}}.
\end{aligned}
\label{estim2uc}
\end{equation}
From (\ref{finh2}), (\ref{estim3u}), (\ref{estim4u}), (\ref{estim1uc}) and (\ref{estim2uc}) we obtain finally from the Plancherel Theorem that there exist $C,\kappa>0$ such that for any $t>0$ we have:
\begin{equation}
\begin{aligned}
& \|\phi(\lambda^{-1}D)W(t)\left(\begin{array}{c}
q_0\\
 m_0\\
\end{array}
\right)\|_{L^2}\leq C e^{-\kappa \lambda^2 t}\|\phi(\lambda^{-1}D)\left(\begin{array}{c}
q_0\\
 m_0\\
\end{array}
\right)\|_{L^2} 
\end{aligned}
\label{estim4ufin}
\end{equation}
It concludes the proof of the proposition \ref{propocru}. {\hfill $\Box$}
\\
We are now going to prove time decay estimates in Besov spaces for the semi group $W(t)$.
\begin{proposition}
Let $s\in\R$, $r\in[1,+\infty]$ and $s_1>s$, $(q_0,m_0)\in( B^{s}_{2,\infty})^{N+1}$ then it exists $C_{s_1}>0$ such that for all $t>0$ we have:
\begin{equation}
\|W(t)\left(\begin{array}{c}
q_0\\
 m_0\\
\end{array}
\right)\|_{B^{s_1}_{2,r}}\leq\frac{C_{s_1}}{t^{\frac{s_1-s}{2}}}\|\left(\begin{array}{c}
q_0\\
 m_0\\
\end{array}
\right)\|_{B^{s}_{2,r}}.
\label{decay1}
\end{equation}
\label{decay}
\end{proposition}
{\bf Proof:} From proposition \ref{propocru}, we have for $\kappa,C>0$ and any $l\in\mathbb{Z}$:
\begin{equation}
\|\D_l W(t)\left(\begin{array}{c}
q_0\\
 m_0\\
\end{array}
\right)\|_{L^2}\leq Ce^{-\kappa 2^{2l}t}\|\D_l \left(\begin{array}{c}
q_0\\
 m_0\\
\end{array}
\right)\|_{L^2}.
\label{decay1a}
\end{equation}
We deduce that for any $l\in\mathbb{Z}$, we have:
\begin{equation}
2^{l s_1}\|\D_l W(t)\left(\begin{array}{c}
q_0\\
 m_0\\
\end{array}
\right)\|_{L^2}\leq \frac{C}{t^{\frac{s_1-s}{2}}} (2^{2l}t)^{\frac{s_1-s}{2}}e^{-\kappa 2^{2l}t} 2^{ls} \|\D_l\left(\begin{array}{c}
q_0\\
 m_0\\
\end{array}
\right)\|_{L^2}.
\label{decay2a}
\end{equation}
We use now the fact that $x^{\frac{s_1-s}{2}} e^{-\kappa x}$ is bounded in $L^\infty(\R^+)$ to deduce that it exists $C_{s_1}>0$ such that:
\begin{equation}
\|2^{l s_1}\|\D_l W(t)\left(\begin{array}{c}
q_0\\
 m_0\\
\end{array}
\right)\|_{L^2}\|_{l^r}\leq \frac{C_{s_1}}{t^{\frac{s_1-s}{2}}} \|2^{ls}\|\D_l \left(\begin{array}{c}
q_0\\
 m_0\\
\end{array}
\right)\|_{l^r}.
\label{decay3a}
\end{equation}
It concludes the proof of the proposition \ref{decay}.  {\hfill $\Box$}
 \subsection{Proof of the Theorem \ref{theo2}}
We shall use a contracting mapping argument to prove the Theorem \ref{theo2} \footnote{For the moment we only consider the case $0<\kappa^2<\mu^2$ and to simplify the notation we assume that $\bar{\rho}=1$.} and we consider the following map $\psi_1$ defined as follows with $q=\rho-1$:
 \begin{equation}
\begin{aligned}
&\psi_1(q,\rho v_1,\rho v_2)=W_{c_1,\mu,\sqrt{\mu^2-\kappa^2},P'(1)}(t)\left(\begin{array}{c}
q_0\\
\rho v_1(0,\cdot)\\
\end{array}
\right)\\
&\hspace{3cm}+\int_{0}^{t}W_{c_1,\mu,\sqrt{\mu^2-\kappa^2},P'(1)}(t-s)\left(\begin{array}{c}
 0 \\
F(\rho,v_1,v_2)\\
\end{array}
\right)(s)\ ds.\\[1mm]
&\psi_2(q,\rho v_1,\rho v_2)=W_{c_2,\mu,-\sqrt{\mu^2-\kappa^2},P'(1)}(t)\left(\begin{array}{c}
q_0\\
\rho v_2(0,\cdot)\\
\end{array}
\right)\\
&\hspace{3cm}+\int_{0}^{t}W_{c_2,\mu,-\sqrt{\mu^2-\kappa^2},P'(1)}(t-s)\left(\begin{array}{c}
 0 \\
F(\rho,v_1,v_2)\\
\end{array}
\right)(s)\ ds.
\end{aligned}
\label{Duhamel1}
\end{equation}
with 
$$
\begin{aligned}
&F(\rho,v_1,v_2)=-\frac{1}{2}({\rm div}(\rho v_1\otimes v_2)+{\rm div}(\rho v_2\otimes v_1))-(P'(\rho)-P'(1))\n\rho.
\end{aligned}
$$
We define finally $\psi_3$ as follows:
\begin{equation}
\begin{aligned}
&\psi_3(q,\rho v_1,\rho v_2)=\left(\begin{array}{c}
\psi_1(q,\rho v_1,\rho v_2)\\
(\psi_2)_2(q,\rho v_1,\rho v_2)\\
\end{array}
\right).
\end{aligned}
\end{equation}
Let us prove now that $\psi_3$ is a map from $X_{\N}$ in itself with $s_1\in(\frac{3}{4},1)$:
$$
\begin{aligned}
&E_\N=\big(\widetilde{L}^\infty(\R^+,B^{\N-1}_{2,\infty}\cap B^{\N}_{2,\infty})\cap \widetilde{L}^1(\R^+,B^{\N+1}_{2,\infty}\cap B^{\N+2}_{2,\infty})\big)\\
&\hspace{6cm}\times \big(\widetilde{L}^\infty(\R^+,B^{\N-1}_{2,\infty})\cap \widetilde{L}^1(\R^+,B^{\N+1}_{2,\infty})\big)^{2N}\\
&\|(q,\rho v_1,\rho v_2)\|_{W_\N}=\sup_{t\in\R^+}  t^{\frac{s_1}{2}}\|(q(t,\cdot),\rho v_1(t,\cdot),\rho v_2(t,\cdot)\|_{B^{\N-1+s_1}_{2,\infty}}\\
&\|(q,\rho v_1,\rho v_2)\|_{X_\N}=\|(q,\rho v_1,\rho v_2)\|_{E_\N\cap W_\N}+\|q\|_{L^\infty(\R^+,L^\infty(\R^N))}.
\end{aligned}
$$
The space $X_\N$ in which we work is more complicated as in \cite{fDD}, indeed we need decay estimate in time in Besov space on the solution in order to control the $L^\infty$ norm of $q$. In \cite{fDD}, the control of the $L^\infty$ norm of $q$ is a direct consequence of Besov embedding $B^{\N}_{2,1}\h L^\infty$ since the third index of the Besov spaces are $1$. From the proposition (\ref{flinear3}), we deduce that for $C>0$ large enough:
\begin{equation}
\begin{aligned}
&\|\psi_3(q,\rho v_1,\rho v_2)\|_{E_\N}\leq C \big(\|q_0\|_{B^{\N-1}_{2,\infty}\times B^{\N}_{2,\infty}}+\|\rho v_1(0,\cdot)\|_{B^{\N-1}_{2,\infty}}+\|\rho v_2(0,\cdot)\|_{B^{\N-1}_{2,\infty}}\\
&\hspace{8cm}+\|F(\rho,v_1,v_2)\|_{\widetilde{L}^1(\R^+,B^{\N-1}_{2,\infty})}\big).
\end{aligned}
\label{EN}
\end{equation}
Next using classical paraproduct law and composition theorems (see \cite{BCD}), we get for $C_1>0$ large enough and a continuous function $C$ using the fact that
${\rm div}(\rho v_1\otimes v_2)=\frac{1}{\rho}\rho v_1\n(\rho v_2)+\frac{1}{\rho}{\rm div}(\rho v_1)\rho v_2+\n(\frac{1}{\rho}-1)\cdot (\rho v_1)\;\rho v_2$:
\begin{equation}
\begin{aligned}
&\|{\rm div}(\rho v_1\otimes v_2)\|_{\widetilde{L}^1(\R^+,B^{\N-1}_{2,\infty})}\leq C_1\big(\|(q,\rho v_1,\rho v_2)\|_{E_\N}^2(1+\|\frac{1}{\rho}-1\|_{L^\infty(\R^+,L^\infty(\R^N))\cap\widetilde{L}^\infty(\R^+,B^{\N}_{2,\infty}))})\\
&\hspace{3cm}+\|(q,\rho v_1,\rho v_2)\|_{E_\N}^2\|\n(\frac{1}{\rho}-1)\|_{\widetilde{L}^{3}(B^{\N-\frac{1}{3}}_{2,\infty})}\big),\\
&\leq C_1\|(q,\rho v_1,\rho v_2)\|_{E_\N}^2(1
+C(\|(\rho,\frac{1}{\rho})\|_{L^\infty(\R^+,L^\infty(\R^N))})\|(q,\rho v_1,\rho v_2)\|_{E_\N}\big)\\[2mm]
&\|(P'(\rho)-P'(1))\n q\|_{\widetilde{L}^1(\R^+,B^{\N-1}_{2,\infty})}\leq C_1 \|\n q\|_{\widetilde{L}^{\frac{4}{3}}(\R^+,B^{\N-\frac{1}{2}}_{2,\infty})} \|P'(1+q)-P'(1)\|_{\widetilde{L}^{4}(\R^+,B^{\N-\frac{1}{2}}_{2,\infty})}\\
&\hspace{4cm}\leq \|(q,\rho v_1,\rho v_2)\|_{E_\N}^2 C(\|q\|_{L^\infty(\R^+,L^\infty(\R^N))}).
\end{aligned}
\label{r1}
\end{equation}
Combining (\ref{EN}), (\ref{r1}), interpolation and composition theorems, we obtain for $C_1>0$ large enough and a continuous function $C$:
\begin{equation}
\begin{aligned}
&\|\psi_3(q,\rho v_1,\rho v_2)\|_{E_\N}\leq C_1 \big(\|q_0\|_{B^{\N-1}_{2,\infty}\times B^{\N}_{2,\infty}}+\|\rho v_1(0,\cdot)\|_{B^{\N-1}_{2,\infty}}+\|\rho v_2(0,\cdot)\|_{B^{\N-1}_{2,\infty}}\\
&\|(q,\rho v_1,\rho v_2)\|^2_{X_\N}(1+C(\|(\rho,\frac{1}{\rho})\|_{L^\infty(\R^+,L^\infty(\R^N))}(1+\|(q,\rho v_1,\rho v_2)\|_{E_\N})\big).
\end{aligned}
\label{ENv}
\end{equation}
It remains now to estimate the $\|\psi_3(q,\rho v_1,\rho v_2)\|_{W^{\N}}$ norm, using the definition of $\psi_1$ and $\psi_2$ (see (\ref{Duhamel1})) we have for $t>0$:
\begin{equation}
\begin{aligned}
&\|\psi_1(q,\rho v_1,\rho v_2)(t,\cdot)\|_{B^{\N-1+s_1}_{2,\infty}}\leq \|W_{c_1,\mu,\sqrt{\mu^2-\kappa^2},P'(1)}(t)\left(\begin{array}{c}
q_0\\
\rho v_1(0,\cdot)\\
\end{array}
\right)\|_{B^{\N-1+s_1}_{2,\infty}}\\
&\hspace{3cm}+\int_{0}^{t}\|W_{c_1,\mu,\sqrt{\mu^2-\kappa^2},P'(1)}(t-s)\left(\begin{array}{c}
 0 \\
F(\rho,v_1,v_2)(s)\\
\end{array}
\right)\|_{B^{\N-1+s_1}_{2,\infty}}\ ds.\\[2mm]
&\|\psi_2(q,\rho v_1,\rho v_2)(t,\cdot)\|_{B^{\N-1+s_1}_{2,\infty}}\leq \|W_{c_1,\mu,-\sqrt{\mu^2-\kappa^2},P'(1)}(t)\left(\begin{array}{c}
q_0\\
\rho v_2(0,\cdot)\\
\end{array}
\right)\|_{B^{\N-1+s_1}_{2,\infty}}\\
&\hspace{3cm}+\int_{0}^{t}\|W_{c_1,\mu,-\sqrt{\mu^2-\kappa^2},P'(1)}(t-s)\left(\begin{array}{c}
 0 \\
F(\rho,v_1,v_2)(s)\\
\end{array}
\right)\|_{B^{\N-1+s_1}_{2,\infty}}\ ds.
\end{aligned}
\label{Duhamel2}
\end{equation}
From proposition \ref{decay}, it yields that for $C>0$ large enough we have:
\begin{equation}
\begin{aligned}
&\|\psi_3(q,\rho v_1,\rho v_2)(t,\cdot)\|_{B^{\N-1+s_1}_{2,\infty}}\leq\frac{C}{t^{\frac{s_1}{2}}}\|(q_0,\rho v_1(0,\cdot),\rho v_2(0,\cdot))\|_{B^{\N-1}_{2,\infty}}\\
&+\int^t_0\frac{C}{(t-s)^{1-\frac{s_1}{2}}}\|F(\rho,v_1,v_2)(s,\cdot)\|_{B^{\N-3+2s_1}_{2,\infty}} ds
\end{aligned}
\label{A1}
\end{equation}
We have now using classical paraproduct laws for $s>0$ and $C$ a continuous function:
\begin{equation}
\begin{aligned}
&\|{\rm div}(\rho v_1\otimes v_2)(s,\cdot)\|_{B^{\N-3+2s_1}_{2,\infty}}\leq\|(\frac{1}{\rho}-1) \rho v_1\otimes \rho v_2)(s,\cdot)\|_{B^{\N-2+2s_1}_{2,\infty}}\\
&\hspace{8cm}+\|(\rho v_1\otimes \rho v_2)(s,\cdot)\|_{B^{\N-2+2s_1}_{2,\infty}}\\
&\leq C_1  C(\|\frac{1}{\rho}(s,\cdot)\|_{L^\infty(\R^N)})(1+\|q(s,\cdot)\|_{B^{\N}_{2,\infty}\cap L^\infty})\|\rho v_1(s,\cdot)\|_{B^{\N-1+s_1}_{2,\infty}}\|\rho v_2(s,\cdot)\|_{B^{\N-1+s_1}_{2,\infty}}
\end{aligned}
\label{parafi1}
\end{equation}
Similarly we have for $C>0$ large enough and $s_1\in(\frac{3}{4},1)$:
\begin{equation}
\begin{aligned}
&\|(P'(\rho)-P'(1))\n q(s,\cdot)\|_{B^{\N-3+2s_1}_{2,\infty}}
\leq C \|\n q(s,\cdot)\|_{B^{\N-2+s_1}_{2,\infty}}\|(P'(1+q)-P'(1))(s,\cdot)\|_{B^{\N-1+s_1}_{2,\infty}}\\
&\leq C(\|q(s,\cdot)\|_{L^\infty}) \|q(s,\cdot)\|^2_{B^{\N-1+s_1}_{2,\infty}}.
\end{aligned}
\label{parafi2}
\end{equation}
From (\ref{A1}), (\ref{parafi1}) and (\ref{parafi2}), we get for $C_1>0$ large enough and $C$ a continuous function:
$$
\begin{aligned}
&\|\psi_3(q,\rho v_1,\rho v_2)(t,\cdot)\|_{B^{\N-1+s_1}_{2,\infty}}\leq \frac{C_1}{t^{\frac{s_1}{2}}}\|(q_0,\rho v_1(0,\cdot),\rho v_2(0,\cdot))\|_{B^{\N-1}_{2,\infty}}\\
&+\int^t_0\frac{1}{(t-s)^{1-\frac{s_1}{2}}}\big(   C(\|\frac{1}{\rho}(s,\cdot)\|_{L^\infty(\R^N)})(1+\|q(s,\cdot)\|_{B^{\N}_{2,\infty}\cap L^\infty})\|\rho v_1(s,\cdot)\|_{B^{\N-1+s_1}_{2,\infty}}\\
&\times\|\rho v_2(s,\cdot)\|_{B^{\N-1+s_1}_{2,\infty}}+C(\|q(s,\cdot)\|_{L^\infty}) \|q(s,\cdot)\|^2_{B^{\N-1+s_1}_{2,\infty}}\big) ds.
\end{aligned}
$$
We get finally for a continuous function $C_2,C_3$:
\begin{equation}
\begin{aligned}
&\|\psi_3(q,\rho v_1,\rho v_2)(t,\cdot)\|_{B^{\N-1+s_1}_{2,\infty}} \leq \frac{C_1}{t^{\frac{s_1}{2}}}\|(q_0,\rho v_1(0,\cdot),\rho v_2(0,\cdot))\|_{B^{\N-1}_{2,\infty}}\\
&+\|(q,\rho v_1,\rho v_2)\|^2_{W^{\N}} \int^t_0\frac{1}{(t-s)^{1-\frac{s_1}{2}}}\frac{1}{s^{s_1}} C_2(\|(\frac{1}{\rho}(s,\cdot),\rho(s,\cdot))\|_{L^\infty(\R^N)})(1+\|q(s,\cdot)\|_{B^{\N}_{2,\infty}\cap L^\infty}) ds\\
&\leq \frac{C_1}{t^{\frac{s_1}{2}}}\|(q_0,\rho v_1(0,\cdot),\rho v_2(0,\cdot))\|_{B^{\N-1}_{2,\infty}}\\
&+\|(q,\rho v_1,\rho v_2)\|^2_{W^{\N}}C_3(\|(\rho,\frac{1}{\rho})\|_{L^\infty(\R^+,L^\infty)})(1+\|(q,\rho v_1,\rho v_2)\|_{X_{\N}})  \int^t_0\frac{1}{(t-s)^{1-\frac{s_1}{2}}}\frac{1}{s^{s_1}} ds
\end{aligned}
\label{A2}
\end{equation}
We have now since $s_1\in(\frac{3}{4},1)$ and for $C>0$ large enough
$$ \int^t_0\frac{1}{(t-s)^{1-\frac{s_1}{2}}}\frac{1}{s^{s_1}} ds=\frac{1}{t^{\frac{s_1}{2}}}\int^1_0\frac{1}{(1-u)^{1-\frac{s_1}{2}}}\frac{1}{u^{s_1}} du\leq \frac{C}{t^{\frac{s_1}{2}}}. 
 $$
 It comes that for $C>0$ large enough and a continuous function $C_1$ we have:
 \begin{equation}
 \begin{aligned}
 &\|\psi_3(q,\rho v_1,\rho v_2)\|_{W_{\N}}\leq C\big(\|(q_0,\rho v_1(0,\cdot),\rho v_2(0,\cdot))\|_{B^{\N-1}_{2,\infty}}\\
 &\hspace{2cm}+\|(q,\rho v_1,\rho v_2)\|^2_{W_{\N}}C_1(\|(\rho,\frac{1}{\rho})\|_{L^\infty(\R^+,L^\infty)})(1+\|(q,\rho v_1,\rho v_2)\|_{X_{\N}}) \big).
 \end{aligned}
 \label{retu1}
 \end{equation}
 It remains now to estimate the $L^\infty(\R^+,L^\infty(\R^N))$ norm on $[\psi_3(q,\rho v_1,\rho v_2)]_1$. From (\ref{Duhamel1}) and using the fact that $B^{\N}_{2,1}(\R^N)$ is embedded in $L^\infty(\R^N)$, we have for $C>0$ large enough:
\begin{equation}
\begin{aligned}
&\|[\psi_3(q,\rho v_1,\rho v_2)]_1(t,\cdot)\|_{L^\infty(\R^+,L^\infty(\R^N))}\leq \|[W_{c_1,\mu,\sqrt{\mu^2-\kappa^2},P'(1)}(t)\left(\begin{array}{c}
q_0\\
\rho v_1(0,\cdot)\\
\end{array}
\right)]_1\|_{L^\infty(\R^N)}  \\
&\hspace{3cm}+\int_{0}^{t}\|[W_{c_1,\mu,\sqrt{\mu^2-\kappa^2},P'(1)}(t-s)\left(\begin{array}{c}
 0 \\
F(\rho,v_1,v_2)\\
\end{array}
\right)]_1(s)\|_{B^{\N}_{2,1}}\ ds.
\end{aligned}
\label{Duhamel1n}
\end{equation}
$W_{c_1,\mu,\sqrt{\mu^2-\kappa^2},P'(1)}\left(\begin{array}{c}
q_0\\
\rho v_1(0,\cdot)\\
\end{array}
\right)$ is the solution of the system (\ref{klinear}), we deduce from the first equation of (\ref{klinear}) and using the maximum principle and the proposition \ref{decay} that for $t>0$ and $C,C_1,C_2,C_3>0$ large enough:
$$
\begin{aligned}
&\|[W_{c_1,\mu,\sqrt{\mu^2-\kappa^2},P'(1)}\left(\begin{array}{c}
q_0\\
\rho v_1(0,\cdot)\\
\end{array}
\right)]_1(t,\cdot)\|_{L^\infty(\R^N)} \leq \|q_0\|_{L^\infty(\R^N)}\\
&\hspace{1cm}+\int^t_0  \|e^{c_1(t-s)\D}{\rm div}
[W_{c_1,\mu,\sqrt{\mu^2-\kappa^2},P'(1)}(s)\left(\begin{array}{c}
q_0\\
\rho v_1(0,\cdot)\\
\end{array}
\right)]_2 \|_{B^{\N}_{2,1}} ds\\[2mm]
& \leq \|q_0\|_{L^\infty(\R^N)}
+\int^t_0 \frac{C}{(t-s)^{1-\frac{s_1}{4}}}\|
[W_{c_1,\mu,\sqrt{\mu^2-\kappa^2},P'(1)}(s)\left(\begin{array}{c}
q_0\\
\rho v_1(0,\cdot)\\
\end{array}
\right)]_2 \|_{B^{\N-1+\frac{s_1}{2}}_{2,1}} ds
\end{aligned}
$$
And it yields:
\begin{equation}
\begin{aligned}
&\|[W_{c_1,\mu,\sqrt{\mu^2-\kappa^2},P'(1)}\left(\begin{array}{c}
q_0\\
\rho v_1(0,\cdot)\\
\end{array}
\right)]_1(t,\cdot)\|_{L^\infty(\R^N)} \leq \|q_0\|_{L^\infty(\R^N)}
\\
&+\int^t_0 \frac{C_1}{(t-s)^{1-\frac{s_1}{4}}}\|
[W_{c_1,\mu,\sqrt{\mu^2-\kappa^2},P'(1)}(s)\left(\begin{array}{c}
q_0\\
\rho v_1(0,\cdot)\\
\end{array}
\right)]_2 \|_{B^{\N-1}_{2,\infty}}^{\frac{1}{2}}\\
&\hspace{4cm}\times \|
[W_{c_1,\mu,\sqrt{\mu^2-\kappa^2},P'(1)}(s)\left(\begin{array}{c}
q_0\\
\rho v_1(0,\cdot)\\
\end{array}
\right)]_2 \|_{B^{\N-1+s_1}_{2,\infty}}^{\frac{1}{2}} ds\\[2mm]
& \leq \|q_0\|_{L^\infty(\R^N)}
+\int^t_0 \frac{C_2}{(t-s)^{1-\frac{s_1}{4}}}\frac{1}{s^{\frac{s_1}{4}}}\|\left(\begin{array}{c}
q_0\\
\rho v_1(0,\cdot)\\
\end{array}
\right)\|_{B^{\N-1}_{2,\infty}} ds\\
& \leq \|q_0\|_{L^\infty(\R^N)}
+C_3 \|\left(\begin{array}{c}
q_0\\
\rho v_1(0,\cdot)\\
\end{array}
\right) \|_{B^{\N-1}_{2,\infty}} .
\end{aligned}
\label{3.100}
\end{equation}
The third inequality corresponds to an interpolation inequality in Besov spaces. In a similar way, we have for $C,C_1,C_2,C_3>0$ large enough and using proposition \ref{decay} and interpolation in Besov spaces (see \cite{BCD}):
\begin{equation}
\begin{aligned}
&\int_0^{t}\|[W_{c_1,\mu,\sqrt{\mu^2-\kappa^2},P'(1)}(t-s)\left(\begin{array}{c}
 0 \\
F(\rho,v_1,v_2)\\
\end{array}
\right)]_1(s)\|_{B^{\N}_{2,1}}\ ds\\
&\leq \int^t_0 \frac{C}{(t-s)^{\frac{3}{2}-\frac{3 s_1}{4}}}\|
\left(\begin{array}{c}
0\\
F(\rho,v_1,v_2)\\
\end{array}
\right) \|_{B^{\N-3+\frac{3}{2} s_1}_{2,1}} ds\\[2mm]
& \leq \int^t_0 \frac{C_1}{(t-s)^{\frac{3}{2}-\frac{3 s_1}{4}}}\big(\|F_2(\rho,v_1,v_2)(s,\cdot)\|^{\frac{2}{3}}_{B^{\N-3+2s_1}_{2,\infty}}\|F_2(\rho,v_1,v_2)(s,\cdot)\|^{\frac{1}{3}}_{B^{\N-3+\frac{s_1}{2}}_{2,\infty}}\\
&+\|F_3(\rho)(s,\cdot)\|^{\frac{1}{4}}_{B^{\N-3+2 s_1}_{2,\infty}}\|F_3(\rho)(s,\cdot)\|^{\frac{3}{4}}_{B^{\N-3+\frac{4 s_1}{3}}_{2,\infty}}\big)
ds,
\end{aligned}
\label{Duhamel1n2}
\end{equation}
with $F_2(\rho,v_1,v_2)=\frac{1}{2}\big({\rm div}(\rho v_1\otimes v_2)+{\rm div}(\rho v_2\otimes v_1)\big)$ and $F_3(\rho)=(P'(\rho)-P'(1))\n\rho$.
Using classical paraproduct law, composition and interpolation theorems (see \cite{BCD}), we have now for $s>0$, a function $C$ continuous and $C_1>0$ large enough:
\begin{equation}\begin{aligned}
&\|{\rm div}(\rho v_1\otimes v_2)(s)\|_{B^{\N-3+\frac{s_1}{2}}_{2,\infty}}\leq\|(\frac{1}{\rho}-1) \rho v_1\otimes \rho v_2)(s,\cdot)\|_{B^{\N-2+\frac{s_1}{2}}_{2,\infty}}\\
&\hspace{8cm}+\|(\rho v_1\otimes \rho v_2)(s,\cdot)\|_{B^{\N-2+\frac{s_1}{2}}_{2,\infty}}\\
&\leq C_1  C(\|\frac{1}{\rho}(s,\cdot)\|_{L^\infty(\R^N)})(1+\|q(s,\cdot)\|_{B^{\N}_{2,\infty}\cap L^\infty})\|\rho v_1(s,\cdot)\|_{B^{\N-1}_{2,\infty}}\|\rho v_2(s,\cdot)\|_{B^{\N-1+\frac{s_1}{2}}_{2,\infty}}\\
&\leq C_1  C(\|\frac{1}{\rho}(s,\cdot)\|_{L^\infty(\R^N)})(1+\|q(s,\cdot)\|_{B^{\N}_{2,\infty}\cap L^\infty})\|\rho v_1(s,\cdot)\|_{B^{\N-1}_{2,\infty}}\\
&\hspace{6cm}\times\|\rho v_2(s,\cdot)\|_{B^{\N-1}_{2,\infty}}^{\frac{1}{2}}
\|\rho v_2(s,\cdot)\|_{B^{\N-1+s_1}_{2,\infty}}^{\frac{1}{2}}.
\end{aligned}
\label{parafi1Ta}
\end{equation}
Similarly we have using composition, interpolation theorems and paraproduct laws (we can since $s_1>\frac{3}{4}$), it exists a continuous function $C$ and $C_1>0$ large enough such that: 
\begin{equation}\begin{aligned}
\|(P'(1+q)-P'(1))&\n q(s,\cdot)\|_{B^{\N-3+\frac{4 s_1}{3}}_{2,\infty}}\\
&\leq \|\n q(s,\cdot)\|_{B^{\N-2+\frac{2s_1}{3}}_{2,\infty}}\|\big(P'(1+q)-P'(1)\big)(s,\cdot)\|_{B^{\N-1+\frac{2s_1}{3}}_{2,\infty}}\\
& \leq C(\|\rho(s,\cdot)\|_{L^\infty(\R^N)}) \| q(s,\cdot)\|^2_{B^{\N-1+\frac{2s_1}{3}}_{2,\infty}}\\
& \leq C_1 C(\|\rho(s,\cdot)\|_{L^\infty(\R^N)}) \| q(s,\cdot)\|^{\frac{4}{3}}_{B^{\N-1+s_1}_{2,\infty}}\| q(s,\cdot)\|^{\frac{2}{3}}_{B^{\N-1}_{2,\infty}}.
\end{aligned}
\label{parafi2Ta}
\end{equation}
It yields using (\ref{parafi1}), (\ref{parafi2}), (\ref{parafi1Ta}) and (\ref{parafi2Ta}) that for $C$ a continuous function and $C_1>0$ sufficiently large:
\begin{equation}
\begin{aligned}
&\|F_2(\rho,v_1,v_2)(s,\cdot)\|^{\frac{2}{3}}_{B^{\N-3+2s_1}_{2,\infty}}\|F_2(\rho,v_1,v_2)(s,\cdot)\|^{\frac{1}{3}}_{B^{\N-3+\frac{s_1}{2}}_{2,\infty}}\\
&\hspace{3cm}+\|F_3(\rho)(s,\cdot)\|^{\frac{1}{4}}_{B^{\N-3+2 s_1}_{2,\infty}}\|F_3(\rho)(s,\cdot)\|^{\frac{3}{4}}_{B^{\N-3+\frac{4 s_1}{3}}_{2,\infty}}\\
&\leq \biggl(C(\|(\frac{1}{\rho},\rho)(s,\cdot)\|_{L^\infty(\R^N)})(1+\|q(s,\cdot)\|_{B^{\N}_{2,\infty}\cap L^\infty})\big(\|\rho v_1(s,\cdot)\|_{B^{\N-1+s_1}_{2,\infty}}\|\rho v_2(s,\cdot)\|_{B^{\N-1+s_1}_{2,\infty}}
\big)\biggl)^{\frac{2}{3}}\times\\
&\biggl(
 C(\|(\frac{1}{\rho},\rho)(s,\cdot)\|_{L^\infty(\R^N)})(1+\|q(s,\cdot)\|_{B^{\N}_{2,\infty}\cap L^\infty})\big( \|\rho v_1(s,\cdot)\|_{B^{\N-1}_{2,\infty}}\|\rho v_2(s,\cdot)\|_{B^{\N-1}_{2,\infty}}^{\frac{1}{2}}
\|\rho v_2(s,\cdot)\|_{B^{\N-1+s_1}_{2,\infty}}^{\frac{1}{2}}\big)\biggl)^{\frac{1}{3}}\\
&+\biggl(C(\|(\frac{1}{\rho},\rho)(s,\cdot)\|_{L^\infty(\R^N)})(1+\|q(s,\cdot)\|_{B^{\N}_{2,\infty}\cap L^\infty}) \|q(s,\cdot)\|^2_{B^{\N-1+s_1}_{2,\infty}}\biggl)^{\frac{1}{4}}\\
&\hspace{3cm}\times \biggl( C_1 C(\|\rho(s,\cdot)\|_{L^\infty(\R^N)}) \| q(s,\cdot)\|^{\frac{4}{3}}_{B^{\N-1+s_1}_{2,\infty}}\| q(s,\cdot)\|^{\frac{2}{3}}_{B^{\N-1}_{2,\infty}}\biggl)^{\frac{3}{4}},\\
&\leq   C(\|(\frac{1}{\rho},\rho)(s,\cdot)\|_{L^\infty(\R^N)})(1+\|q\|_{E_{\N}\cap L^\infty(\R^+,L^\infty(\R^N))})\biggl( \frac{1}{s^{\frac{3 s_1}{4}}}\|(q,\rho v_1,\rho v_2)\|^{\frac{3}{2}}_{W_{\N}}
\\
&\hspace{3cm}\times \|(q,\rho v_1,\rho v_2)\|^{\frac{1}{2}}_{E_{\N}}+\frac{1}{s^{\frac{3s_1}{4}}}\|(q,\rho v_1,\rho v_2)\|_{W_{\N}}^{\frac{3}{2}}\|(q,\rho v_1,\rho v_2)\|^{\frac{1}{2}}_{E_{\N}}\biggl),
\\
&\leq   \frac{1}{s^{\frac{3s_1}{4}}} C(\|(\frac{1}{\rho},\rho)(s,\cdot)\|_{L^\infty(\R^N)})(1+\|q\|_{E_{\N}\cap L^\infty(\R^+,L^\infty(\R^N))})\|(q,\rho v_1,\rho v_2)\|^{2}_{X_{\N}}.
\end{aligned}
\label{mtech}
\end{equation}
Plugging (\ref{mtech}) in (\ref{Duhamel1n2}), we obtain for a continuous fonction $C$ since $s_1>\frac{3}{4}$:
\begin{equation}
\begin{aligned}
&\int_0^{t}\|[W_{c_1,\mu,\sqrt{\mu^2-\kappa^2},P'(1)}(t-s)\left(\begin{array}{c}
 0 \\
F(\rho,v_1,v_2)\\
\end{array}
\right)]_1(s)\|_{B^{\N}_{2,1}}\ ds\\
& \leq \frac{C_1}{\sqrt{t}}C(\|(\frac{1}{\rho},\rho)(s,\cdot)\|_{L^\infty(\R^N)})(1+\|(q,\rho v_1,\rho v_2)\|_{X_{\N}})\|(q,\rho v_1,\rho v_2)\|^{2}_{X_{\N}}.
\end{aligned}
\label{Duhamel1n3}
\end{equation}
Combining (\ref{Duhamel1n}), (\ref{3.100}) and (\ref{Duhamel1n3}), we have for $C_1>0$ sufficiently large and $C$ a continuous function we have for any $t>0$:
\begin{equation}
\begin{aligned}
&\|[\psi_3(q,\rho v_1,\rho v_2)]_1(t,\cdot)\|_{L^\infty(\R^+,L^\infty(\R^N))}  \leq \|q_0\|_{L^\infty(\R^N)}
+C_1 \|\left(\begin{array}{c}
q_0\\
\rho v_1(0,\cdot)\\
\end{array}
\right) \|_{B^{\N-1}_{2,\infty}} \\
&+\frac{1}{\sqrt{t}}C(\|(\frac{1}{\rho},\rho)(s,\cdot)\|_{L^\infty(\R^N)})(1+\|(q,\rho v_1,\rho v_2)\|_{X_{\N}})\|(q,\rho v_1,\rho v_2)\|_{W_{\N}} \|(q,\rho v_1,\rho v_2)\|_{E_{\N}}.
\end{aligned}
\label{3.101}
\end{equation}
This estimate proves that we need to distinguish the short time and the long time, indeed (\ref{3.101}) is interesting only for the long time.\\
We are going now to apply a fixed point theorem and we define the map $\psi_4$ as follows for any $t>0$:
$$
\begin{cases}
\begin{aligned}
&\psi_4(q,\rho v_1,\rho v_2)(t,\cdot)=\psi(q,\rho v_1,\rho v_2)(t,\cdot)\;\;\;\mbox{if}\;t\in(0,T],\\
&\psi_4(q,\rho v_1,\rho v_2)(t,\cdot)=\psi_3(q,\rho v_1,\rho v_2)(t,\cdot)\;\;\;\mbox{if}\;t>T,
\end{aligned}
\end{cases}
$$
with $\psi$ defined as in the previous section:
$$
\begin{aligned}
&\psi(q,\rho v_1,\rho v_2)=\left(\begin{array}{c}
e^{c_1t\D}q_0\\
e^{\mu t\D}\mathbb{P}(\rho v_1(0,\cdot))+e^{c_2 t\D}\mathbb{Q}(\rho v_1(0,\cdot))\\
e^{\mu t\D}\mathbb{P}(\rho v_2(0,\cdot))+e^{c_1 t\D}\mathbb{Q}(\rho v_2(0,\cdot))\\
\end{array}
\right)\\
&-\int_{0}^{t}\left(\begin{array}{c}
 e^{c_1(t-s)\D}{\rm div}(\rho v_1) \\
 \frac{1}{2} e^{\mu(t-s)\D} \mathbb{P}({\rm div}(\rho v_1\otimes v_2)+{\rm div}(\rho v_2\otimes v_1)) +e^{c_2(t-s)\D} \mathbb{Q}(F(\rho,v_1,v_2))\\
 \frac{1}{2} e^{\mu(t-s)\D} \mathbb{P}({\rm div}(\rho v_1\otimes v_2)+{\rm div}(\rho v_2\otimes v_1)) +e^{c_1(t-s)\D} \mathbb{Q}(F(\rho,v_1,v_2)) \\
\end{array}
\right)\ ds.
\end{aligned}
$$
 and $T>0$ to determine later. We are going now to prove that $\psi$ is a map from $X_{\N,R}\cap E_{R_1,M,T}$ in itself with $s_1\in(\frac{3}{4},1)$ and $0<R,R_1,T<\frac{1}{2}$ sufficiently small that we will define later.
 We define $X_{\N,R}$ as follows:
$$
\begin{aligned}
&X_{\N,R}=\{(q,\rho v_1,\rho v_2)\in X_{\N},\|(q,\rho v_1,\rho v_2)\|_{X_\N}\leq R \}\\
\end{aligned}
$$
In the sequel we will note $X^+_{\N,T}$ and $X^-_{\N,T}$ to define the subsets of $X_{\N}$ where the norms are respectively only considered on the time interval $[T,+\infty[$ and $[0,T]$. It is important to mention that we can work in small time on $(0,T)$ in $E_{R_1,M,T}$ since $(\rho v_1(0,\cdot),\rho v_2(0,\cdot))$ are in $BMO^{-1}(\R^N)$, indeed we know that $B^{\N-1}_{2,\infty}$ is embedded in $BMO^{-1}(\R^N)$.
We have obviously $\|\cdot\|_{X_\N}\leq \|\cdot\|_{X_\N^-}+\|\cdot\|_{X_\N^+}$.\\
Combining (\ref{ENv}),  (\ref{retu1}) and (\ref{3.101}) we have for $C$ a continuous function and $C_1>0$ large enough and for all $(q,\rho v_1,\rho v_2)\in X_{\N,R}\cap E_{R_1,M,T}$:
 \begin{equation}
 \begin{aligned}
 &\|\psi_4 (q,\rho v_1,\rho v_2)\|_{X^+_{\N}}\leq \|\psi_3 (q,\rho v_1,\rho v_2)\|_{X^+_{\N}}
  \leq C_1\big(\|(q_0,\rho v_1(0,\cdot),\rho v_2(0,\cdot))\|_{B^{\N-1}_{2,\infty}}\\
  &+\|q_0\|_{ B^{\N}_{2,\infty}\cap L^\infty}+\frac{1}{\sqrt{T}}\|(q,\rho v_1,\rho v_2)\|^2_{X_{\N}}(1+C(\|\rho,\frac{1}{\rho}\|_{L^\infty(\R^+,L^\infty)}))(1+\|(q,\rho v_1,\rho v_2)\|_{X_{\N}}) \big)\\
&\leq C_1\big(\|(q_0,\rho v_1(0,\cdot),\rho v_2(0,\cdot))\|_{B^{\N-1}_{2,\infty}}+\|q_0\|_{ B^{\N}_{2,\infty}\cap L^\infty}\\
&\hspace{4cm}+\frac{R^2}{\sqrt{T}}(1+\sup_{y\in[0,\max (R+1,\frac{1}{1-R})]} C(y)) (1+R)  \big)
 \end{aligned}
 \label{retu1c}
 \end{equation}
 We take now $R$ such that:
 $$R=2C_1 \big(\|(q_0,\rho v_1(0,\cdot),\rho v_2(0,\cdot))\|_{B^{\N-1}_{2,\infty}}+\|q_0\|_{ B^{\N}_{2,\infty}\cap L^\infty}\big).$$
We assume now that $0<R<\frac{1}{2}$ and then the initial data are sufficiently small such that:
$$
\begin{aligned}
&R^2(1+\sup_{y\in[0,2]} C(y)) (1+R)  \leq\frac{R\sqrt{T} }{2}
\end{aligned}
$$
It gives the following condition:
\begin{equation}
\frac{3}{2}R(1+\sup_{y\in[0,2]} C(y)) \leq\frac{\sqrt{T}}{4}.
\label{petit}
\end{equation}
In a similar way we can estimate $\|\psi_4\|_{X_\N^-\cap E_{R_1,M,T}}$. It suffices to apply exactly the same estimates than in the previous section for the norm $E_{R_1,M,T}$ (with $T$ and $R_1$ small enough). For the norm $X_\N^-$ we have to repeat the same estimates in Besov spaces on a finite time interval $(0,T)$ and to choose $T,R>0$ sufficiently small as previously. We just say few words on the case of the norm $W_{\N}^{-}$. First we have by interpolation since $s_1\in(\frac{3}{4},1)$:
$$
\begin{aligned}
\sup_{t\in[0,T]}t^{\frac{s_1}{2}}\|[\psi_4(q,\rho v_1,\rho v_2)(t,\cdot)]_1\|_{B^{\N-1+s_1}_{2,\infty}}&\leq T^{\frac{s_1}{2}}\|[\psi_4(q,\rho v_1,\rho v_2)(t,\cdot)]_1\|_{\widetilde{L}^\infty([0,T],B^{\N-1}_{2,\infty}\cap B^{\N}_{2,\infty})},\\
&\leq T^{\frac{s_1}{2}}\|\psi_4\|_{E_\N^-}.
\end{aligned}
$$
The norm $E_\N^-$ is classical to estimate (we refer to \cite{BCD} or to the previous estimates) and it is important to mention that we use in a crucial way the fact that the $L^\infty$ norm of $[\psi_4(q,\rho v_1,\rho v_2)]_1$ is bounded by the norm coming from $E_{R_1,M,T}$.
We deduce now that for $C>0$ large enough:
$$
\begin{aligned}
&\sup_{t\in[0,T]}t^{\frac{s_1}{2}}\|\psi_4(q,\rho v_1,\rho v_2)(t,\cdot)\|_{B^{\N-1+s_1}_{2,\infty}}\leq  T^{\frac{s_1}{2}}\|\psi_4\|_{E_\N^-}+C(\|\rho v_1(0,\cdot)\|_{B^{\N-1+s_1}_{2,\infty}}\\
&+\|\rho v_2(0,\cdot)\|_{B^{\N-1+s_1}_{2,\infty}}\\
&\hspace{1cm}+C\int^t_0 \sum_{i=1}^3 \|e^{\beta_i(t-s)\D}[{\rm div}(\rho v_1\otimes v_2+v_2\otimes v_1)(s)+\n P(\rho)(s)] \|_{B^{\N-1+s_1}_{2,\infty}} ds,
\end{aligned}
$$
with $\beta_1\in (c_1,c_2,\mu)$. We have now for $C_1,C_2>0$ large enough, $C$ a continuous function and $t\in[0,T]$:
\begin{equation}
\begin{aligned}
&\int^t_0\|e^{\beta_i(t-s)\D}{\rm div}(\rho v_1\otimes v_2+v_2\otimes v_1)(s) \|_{B^{\N-1+s_1}_{2,\infty}} ds\\
&\leq C_1 \int^t_0\frac{1}{(t-s)^{1-\frac{s_1}{2}}} \|\rho v_1\otimes v_2(s)\|_{B^{\N-2+2s_1}_{2,\infty}}\\
&\leq \int^t_0\frac{1}{(t-s)^{1-\frac{s_1}{2}}}\frac{1}{s^{s_1}}\|(q,\rho v_1,\rho v_2)\|_{W_\N^{-}}^2 C(\|\frac{1}{\rho}(s,\cdot)\|_{L^\infty})(1+\|(q,\rho v_1,\rho v_2)\|_{E_\N^{-}})ds\\
&\leq \frac{C_2}{t^{\frac{s_1}{2}}} \|(q,\rho v_1,\rho v_2)\|_{W_\N^{-}}^2 C(\|\frac{1}{\rho}(s,\cdot)\|_{L^\infty})(1+\|(q,\rho v_1,\rho v_2)\|_{E_\N^{-}}).
\end{aligned}
\label{302}
\end{equation}
Proceeding as previously, we have also  for $C_1,C_2>0$ large enough, $C$ a continuous function and $t\in[0,T]$:
\begin{equation}
\begin{aligned}
\int^t_0\|e^{\beta_i(t-s)\D} \n P(\rho) (s) \|_{B^{\N-1+s_1}_{2,\infty}} ds
&\leq C_2 \int^t_0\frac{1}{(t-s)^{\frac{s_1}{2}}}\|P(1+q)-P(1)\|_{B^\N_{2,\infty}} ds\\
&\leq \frac{C_2 T}{t^{\frac{s_1}{2}}} C(\|\rho(s,\cdot)\|_{L^\infty})\|(q,\rho v_1,\rho v_2)\|_{E_\N^{-}}.
\end{aligned}
\label{303}
\end{equation}
We obtain then the stability for the norm $X_\N^-$ using (\ref{302}) and (\ref{303}). In a similar way, we can prove that $\psi_4$ is contractive. More precisely taking $(q,\rho v_1,\rho v_2)$, $(q_1,\rho_1w_1,\rho_1w_2)$ in $X_{\N,R}\cap E_{R_1,M,T}$ we have for $t>T$:
 \begin{equation}
\begin{aligned}
&\psi_1(q,\rho v_1,\rho v_2)-\psi_1(q_1,\rho_1 w_1,\rho_1 w_2)\\
&\hspace{2cm}=\int_{0}^{t}W_{c_1,\mu,\sqrt{\mu^2-\kappa^2},P'(1)}(t-s)\left(\begin{array}{c}
 0 \\
F(\rho,v_1,v_2)-F(\rho_1,w_1,w_2)\\
\end{array}
\right)(s)\ ds.\\[2mm] 
&\psi_2(q,\rho v_1,\rho v_2)-\psi_2(q_1,\rho_1 w_1,\rho_1 w_2)\\
&\hspace{2cm}=\int_{0}^{t}W_{c_1,\mu,-\sqrt{\mu^2-\kappa^2},P'(1)}(t-s)\left(\begin{array}{c}
 0 \\
F(\rho,v_1,v_2)-F(\rho_1,w_1,w_2)\\\
\end{array}
\right)(s)\ ds.
\end{aligned}
\label{Duhamel1v}
\end{equation}
As previously we show that for $C_1>0$ large enough we have:
\begin{equation}
\begin{aligned}
&\|\psi_4(q,\rho v_1,\rho v_2)-\psi_4(q_1,\rho_1 w_1,\rho_1 w_2)\|_{X^+_{\N}}\\
&\hspace{3cm}\leq C_1 R (1+R)
\|(q,\rho v_1,\rho v_2)-(q_1,\rho_1 w_1,\rho_1 w_2)\|_{X_{\N}}.
\end{aligned}
\end{equation}
Taking again $R$ sufficiently small we deduce that the map $\psi_4$ is contractive (it suffices again to repeat the same process on $[0,T]$ with $T>0$ sufficiently small). We proceed similarly for $X_{\N}^-$ and $E_{R_1,M,T}$.
It concludes the proof of the Theorem \ref{theo2} in the case $0<\kappa^2<\mu^2$. The proof in the case $\kappa^2=\mu^2$ is similar except that $v_1=v_2$, in addition when we study the system (\ref{klinearp}) there is no distinction between high and low frequencies. 
\section*{Acknowledgements}
The author has been partially funded by the ANR project 
INFAMIE ANR-15-CE40-0011. This work was realized during the secondment of the author in the ANGE Inria team.


\begin{thebibliography}{}
\bibitem{fA}
D. M. Anderson, G. B McFadden and A. A. Wheller. Diffuse-interface
methods in
fluid mech. In Annal review of fluid mechanics, Vol. 30, pages 139-165. Annual Reviews, Palo Alto,
CA, 1998.
\bibitem{Aa}
P. Antonelli and S. Spirito, Global existence of finite energy weak solutions of quantum Navier-Stokes equations, \textit{Arch. Rat. Mech. Anal.} 225, no. 3 (2017), 1161-1199.
\bibitem{BCD}
H. Bahouri, J.-Y. Chemin, R. Danchin. Fourier analysis and nonlinear partial differential equations, \textit{Grundlehren der mathematischen Wissenschaften}, \textit{343}, \textit{Springer Verlag}, 2011.
\bibitem{fC}
J. W. Cahn and J. E. Hilliard, Free energy of a nonuniform system, I.
Interfacial free energy, J. Chem. Phys. 28 (1998) 258-267.
\bibitem{CMP}
M. Cannone, Y. Meyer and F. Planchon. Solutions auto-similaires des \'equations de Navier-Stokes. \textit{S\'eminaire sur les \'equations
aux d\'eriv\'ees partielles, 1993-1994}, exp. No12 pp. \'Ecole polytech, palaiseau, 1994.
\bibitem{C}
G.-H Cottet. \'Equations de Navier-Stokes dans le plan avec tourbillon initial mesure. \textit{C. R. Acad.
Sci. Paris S\'er. I Math.} 303(4), 105-108 (1986).
\bibitem{fDD}
R. Danchin and B. Desjardins, Existence of solutions for
compressible fluid models of Korteweg type, \textit{Ann. Inst. Henri Poincar\'e, Analyse Non Lin\'eaire}, 18,97-133 (2001)
\bibitem{fDS}
J. E. Dunn and J. Serrin, On the thermomechanics of interstitial
working, \textit{Arch. Rational Mech. Anal.} 88(2) (1985) 95-133.
\bibitem{FK}
H. Fujita and T. Kato. On the Navier-Stokes initial value problem I, \textit{Arch. Ration. Mech. Anal.},
16, (1964), 269-315.
\bibitem{GG}
I. Gallagher and T. Gallay, Uniqueness for the two-dimensional Navier-Stokes equation
with a measure as initial vorticity, \textit{Math. Ann.}, 332 (2005), 287-327.
\bibitem{GM}
Y. Giga and T. Miyakawa. Navier-Stokes flow in $\R^3$ with measures as initial vorticity and Morrey spaces. \textit{Comm. Partial Differential Equations }14 (1989), 577-618.
\bibitem{GM}
Y. Giga, T. Miyakawa and H. Osada. Two-dimensional Navier-Stokes flow with measures as initial
vorticity. \textit{Arch. Ration. Mech. Anal.} 104(3), 223-250 (1988).
\bibitem{fGP}
M. E. Gurtin, D. Poligone and J. Vinals, Two-phases binary fluids and
immiscible fluids described by an order parameter, \textit{Math. Models
Methods Appl. Sci.} 6(6) (1996) 815--831.
\bibitem{arma}
B. Haspot, Existence of global strong solutions in critical spaces for barotropic viscous fluids,  \textit{Arch. Rational. Mech. Anal},  202, Issue 2 (2011), 427-460.
\bibitem{fH1}
B. Haspot, Existence of solutions for compressible fluid models of Korteweg type, \textit{Annales Math\'ematiques Blaise Pascal} 16, 431-481 (2009).
\bibitem{JMAA}
B. Haspot, Existence of global strong solution for Korteweg system with large infinite energy initial data,  \textit{Journal of Mathematical Analysis and Applications}  438 (2016), pp. 395-443.
\bibitem{Hat} H. Hattori and D. Li. The existence of global solutions to a fluid dynamic model for materials for Korteweg type. \textit{J. Partial Differential Equations}, 9(4): 323-342, 1996.
\bibitem{fJL}
D. Jamet, O. Lebaigue, N. Coutris and J.M. Delhaye, The second
gradient method for the direct numerical simulation of liquid-vapor
flows
with phase change. \textit{J. Comput. Phys}, 169(2): 624--651, (2001).
\bibitem{Ju}
 A. J\"ungel, Global weak solutions to compressible Navier-Stokes equations for quantum fluids, \textit{SIAM journal on mathematical analysis}, 2011, vol. 42, no 3, pp. 1025-1045.
\bibitem{K1}
 T. Kato. The Navier-Stokes equation for an incompressible fluid in $\R^2$ with a measure as the initial
vorticity. \textit{Differ. Integral Eqs.} 7(3Ð4), 949-966 (1994).
\bibitem{fK}
D.J. Korteweg. Sur la forme que prennent les \'equations du
mouvement des fluides si l'on tient compte des forces capillaires
par des variations de densit\'e. \textit{Arch. N\'eer. Sci. Exactes S\'er.
II}, 6 :1-24, 1901.
\bibitem{KT}
H. Koch and D. Tataru. Well-posedness for the Navier-Stokes equations, \textit{Adv. Math.}. 157, 2001, 22-35.
\bibitem{Lemarie}
P.-G. Lemari\'e, Recent developments in the Navier-Stokes problem, \textit{Chapman and Hall/CRC research Notes in Mathematics}, 431 Chapman and Hall/CRC, Boca Raton, FL, 2002
\bibitem{fR}
J. S. Rowlinson, Translation of J.D van der Waals, The thermodynamic
theory of capillarity under the hypothesis of a continuous variation
of density. \textit{J.Statist. Phys.}, 20(2): 197-244, 1979.
\bibitem{fTN}
C. Truedell and W. Noll. The nonlinear field theories of mechanics.
\textit{Springer-Verlag}, Berlin, second edition, 1992.
\bibitem{VW}
J. F. Van der Waals, Thermodynamische Theorie der Kapillarit\"at unter Voraussetzung stetiger Dichte\"anderung, \textit{Phys. Chem.} 13, 657-725 (1894).
\end{thebibliography}
\end{document}